\documentclass[11pt,a4paper]{amsart}
\usepackage{amsmath,amssymb,amsthm,amscd}
\usepackage{enumerate}
\usepackage{graphicx}
\usepackage{psfrag}

\title[Spectral triples]{Spectral triples and aperiodic order}
\thanks{Work supported by the ANR grant {\em SubTile} no. NT09 564112.}
\author{J. Kellendonk \& J. Savinien}
\address{Universit\'e de Lyon,
Universit\'e Lyon 1,
CNRS, UMR 5208 Institut Camille Jordan,
B\^atiment du Doyen Jean Braconnier,
43, blvd du 11 novembre 1918,
F - 69622  Villeurbanne Cedex,
France}                                         

\date{\today}

\newcommand{\N}{\mathbb{N}}
\newcommand{\Z}{\mathbb{Z}}

\newcommand{\R}{\mathbb{R}}
\newcommand{\Rr}{\mathcal{R}}
\newcommand{\C}{\mathbb{C}}

\renewcommand{\P}{\mathbb{P}}
\newcommand{\EE}{\mathbb{E}}

\renewcommand{\L}{\mathcal L}
\newcommand{\E}{\mathcal {E}}

\newcommand{\A}{\mathcal A}
\newcommand{\Ypsilon}{\mathcal Y}
\newcommand{\im}{\mbox{\rm im}\,}
\newcommand{\Tr}{\mbox{\rm Tr}}

\newcommand{\diag}{\mbox{\rm diag}}
\newtheorem{thm}{Theorem}[section]
\newtheorem{lemma}[thm]{Lemma}
\newtheorem{cor}[thm]{Corollary}
\newtheorem{prop}[thm]{Proposition}
\theoremstyle{definition}
\newtheorem{defn}[thm]{Definition}

\theoremstyle{remark}

\numberwithin{equation}{section}

\renewcommand{\H}{\mathcal H}
\newcommand{\Tt}{\mathcal T}
\newcommand{\Pp}{\mathcal P}

\newcommand{\supp}{\mbox{\rm supp}}

\newcommand{\ri}{\underline{r}}
\newcommand{\ro}{\overline{r}}

\renewcommand{\b}{b}

\newcommand{\Aa}{{\mathcal A}}

\newcommand{\Hh}{{\mathcal H}}

\newcommand{\Ll}{{\mathcal L}}

\newcommand{\s}{g}

\newcommand{\NM}{{\mathbb N}}

\newcommand{\PM}{{\mathbb P}}
\newcommand{\RM}{{\mathbb R}}

\newcommand{\punc}{\text{\rm punc}}
\newcommand{\ru}{\overline{\rho}}
\newcommand{\rd}{\underline{\rho}}
\newcommand{\freq}{\text{\rm freq}}

\newcommand{\clr}{C_{\textrm{\tiny LR}}}

\begin{document}

\date{\today}

\maketitle

\begin{abstract}
We construct spectral triples for compact metric spaces $(X,d)$.
%The space is approximated by a graph whose set of vertices is dense in $X$.
This provides us with a new metric $\bar{d}_s$ on $X$.
%the Connes, or spectral distance. %is an extension of the graph metric.
We study its relation with the original metric $d$.
When $X$ is a subshift space, or a discrete tiling space, and $d$ satisfies certain bounds
we advocate
that the property of $\bar{d}_s$ and $d$ to be Lipschitz equivalent is
a characterization of high order. 
For episturmian subshifts, we prove that $\bar{d}_s$ and $d$ are
Lipschitz equivalent if and only if the subshift is repulsive (or power free). 
For Sturmian subshifts this is equivalent to linear recurrence.
For repetitive tilings we show that if their patches have equi-distributed frequencies then the two metrics are  Lipschitz equivalent. Moreover, we
study the zeta-function of the spectral triple and relate its abscissa of convergence to the complexity exponent of the subshift or the tiling.  Finally, we
derive Laplace operators from the spectral triples and compare our construction with 
that of  Pearson and Bellissard.
\end{abstract}

\newpage

\tableofcontents

\newpage

%%%%%%%%%%%%%%%%%%%%%%%%%%%%%%%%%%%%%%%%%%%%%%%%%%%%%%%%%%%%%%%%%%%
\section*{Introduction}
The fundamental notion in non commutative Riemannian geometry is that
of a spectral triple $(\Aa,D,\H)$ for an algebra $\Aa$ \cite{Co94}.
The algebra $\Aa$ acts faithfully on the
Hilbert space $\H$ together with a self-adjoint operator $D$, called
the Dirac operator, which has compact resolvent and bounded commutator
with the elements of $\Aa$. The triple plays the role of the metric for
the possibly virtual space described by $\Aa$.
%, also for non-commutative spaces, i.e.\ virtual spaces defined by non
%commutative algebras. 
In fact it allows to define the distance between two states of the
algebra. Moreover, the commutator
with the Dirac operator leads to the definition of a 
gradient and then of Laplacians via Dirichlet forms. Furthermore the
zeta-function  
$s\mapsto \zeta_D(s):=\Tr(|D|^{-s})$  with its complex analytical
properties ought to carry information about the system the spectral
triple is supposed to describe.  
This is already very interesting in the commutative context, i.e.\ for
spectral triples of algebras of functions on ordinary topological
spaces. Examples of that type 
are fractals, for which the zeta-function plays an important role
\cite{La97,GI03,GI05}.  

The recent work of Pearson and Bellissard \cite{Pea08,PB09} paved the
way to construct 
%on spectral triples for ultra metric Cantor sets  It seems surprisingly hard to construct 
spectral triples for aperiodic tilings, that is, for tiling spaces.
They succeeded in constructing a family of spectral triples for ultra
metric Cantor sets, studied their zeta-function, proposed a family of
Laplacians and initiated the study of the Markov-processes defined by
the Laplacians. This can be directly applied to discrete tiling spaces of
aperiodic tilings. Indeed, the spectral analysis of the family of
Pearson-Bellissard Laplacians was carried out in great detail for
substitution tilings \cite{JS10a}.  
%What does this tell us about aperiodic tilings?

Aperiodic tilings play an important role in the theory of aperiodic
order in that they serve as models for structures which are aperiodic
but still show signs of order. What that really means is still
somewhat under discussion but quite a few concepts have been developed
to characterize order. Patch counting complexity (combinatorial
complexity), repetitivity, and existence and distribution of patch
frequencies are well known concepts to quantify order. They are rather
direct properties of 
tilings. More complicated concepts involve the topological invariants
(finite versus infinite rank of cohomology) or the measure type of the
diffraction measure, for instance whether it is pure point or not. In
view of this it seems natural to ask what spectral triples for
tilings can say about aperiodic order. The present article gives a
first answer to this question. 

The idea is the following. A tiling defines a space $\Xi$. A spectral
triple for the algebra of continuous functions on this space provides
us with a metric on $\Xi$. If we had a canonical way of constructing
such a spectral triple then the properties of this so-called spectral
metric could be understood as properties of the tiling and hence used
to characterize it. Yet we do not have a canonical spectral triple. But
given a metric $d$ on $\Xi$ we have a natural way of constructing a
family of spectral triples on $\Xi$. We thus get a family of new
metrics. %We can take the extremes. 
The infimum over the new metrics yields back the original metric $d$. 
The question of how the supremum over the new metrics relates to the
original one leads to a characterization of $(\Xi,d)$. In particular, it may be
continuous (and hence inducing the same topology) or even equivalent
to $d$, and in the latter case we would say that this is a sign of {\em high order}.
Now there does not seem to be a canonical choice for a metric on tiling spaces.  
We consider therefore a whole family of metrics $d$ on $\Xi$ and
perform the analysis for the family.

\begin{figure}[!h]  
\begin{center}
\psfrag{A}{{\small slope $\theta$}}
\includegraphics[width=9cm]{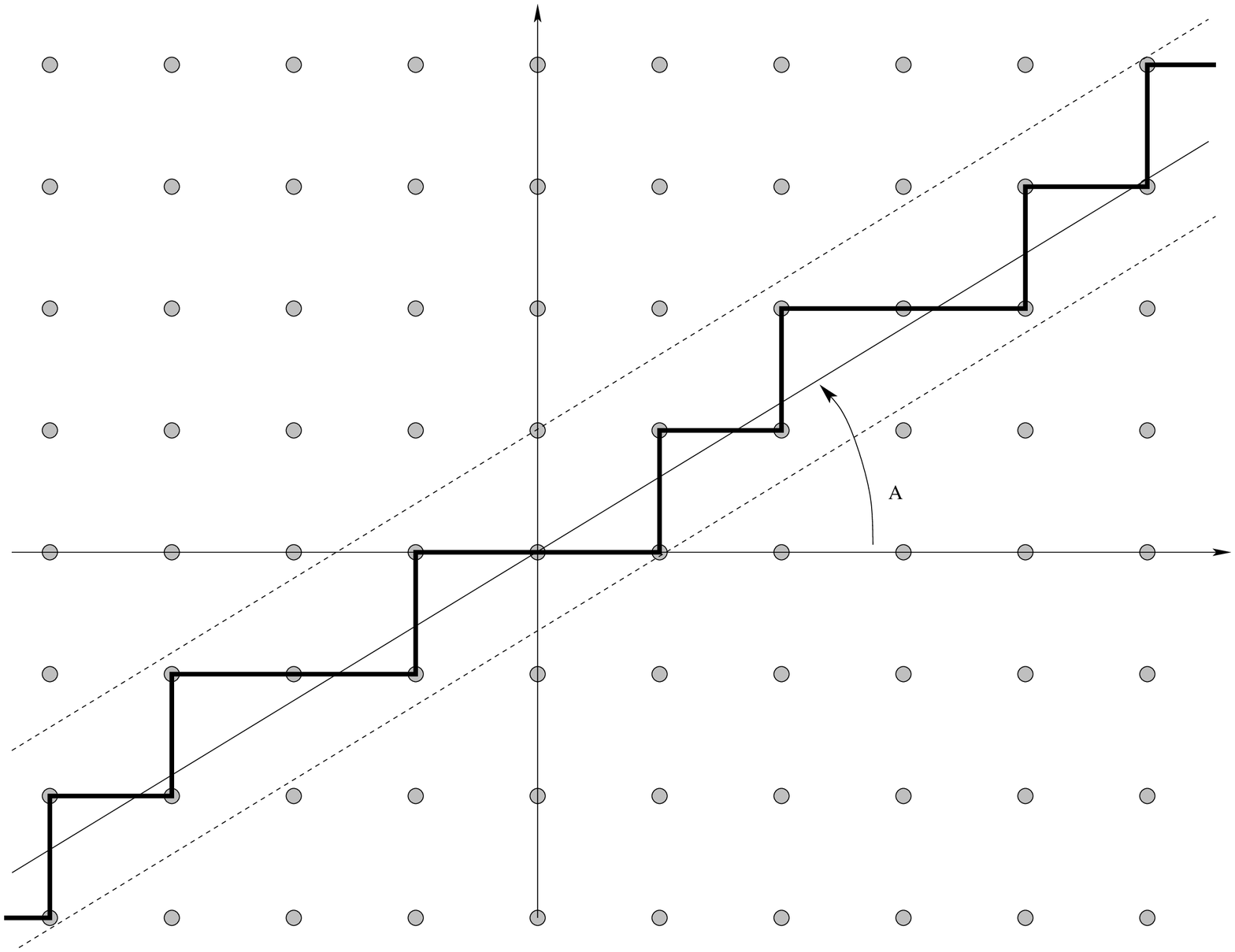}  
\end{center}
\caption{\small{}}
\label{fig-1}
\end{figure}

To give an idea of how this works, we present now the perhaps simplest examples, namely one-sided
canonical cut \& project tilings with dimension and co-dimension equal
to one, see Figure~\ref{fig-1}.
Those tilings correspond to unilateral Sturmian subshifts.

%%%%%%%%%%%%%%%%%%%%%
%%%%%%%%%%%%%%%%%%%%%

The figure suggests that we can approximate the line of irrational
slope $\theta$, say, by a staircase curve. A one-dimensional tiling is
obtained by projecting the lattice points lyoing on the staircase onto the
line, and considering them as the boundaries of tiles. We concentrate on the
right part of it, i.e.\ everything in the right upper quadrant.
We may also read the staircase symbolically:
the horizontal segment as letter $a$, say, and the vertical one as
letter $b$. The associated space $\Xi$ can be described in a
combinatorial way. There is the tree of words
$\Tt$. The vertices of this tree stand for the
different finite words one can find in the sequence. A vertex of level
$n$ represents a word of length $n$. If a word of length $n+1$
extends a word of length $n$ (by a letter at its right end) then we draw
an edge from the shorter to the longer word. These are
precisely the edges of $\Tt$. 
We view them vertically, like in a real tree.
$\Xi$ can be identified with the boundary of $\Tt$, that is, 
it corresponds to the set of infinite paths on the tree, and a path can
be understood as an infinite word. 

We introduce horizontal edges. These are edges between words of the
same length. We draw an un-oriented edge between two distinct words of the same
length if they are extensions (by one letter at the right end) of the same
word. An unoriented edge counts as two oppositely oriented edges.

Given any strictly decreasing sequence
$(\delta_n)_n$ which tends to $0$, we can define a metric on $\Xi$ by
saying that the distance between two (distinct) paths is $\delta_n$ with $n$
chosen to be the length of the longest common prefix of the two.  
On the other hand, we use the sequence to give the horizontal edges a
length, namely we say that an edge $e$ between two words of length $n+1$
has length $l(e):=\delta_n$. 

The final ingredient to define the spectral triple is a collection of
choices, namely we chose for each word one particular extension (by one
letter at its right end). Thus recursively we associate to each word
an infinite extension, i.e.\ an infinite sequence of which it is a
prefix. Call that sequence $\tau(w)\in\Xi$ for the word $w$.

Let $E$ be the set of horizontal edges.  
An edge $e$ has a source $s(e)$ and a range $r(e)$
and these can be identified with vertices, i.e.\ words. 
With $e$, $E$ contains also $\tilde e$, the edge with
opposite orientation.
The Hilbert space of the triple is $\Hh=\ell^2(E)$.
The algebra $\Aa=C(\Xi)$ is represented by
\[\pi_\tau(f) \psi(e)  =  f( \tau(s(e)) ) \psi(e) \]
and the Dirac operator is given by
\[ D\psi(e) = l(e)^{-1}\psi(\tilde{e}).\]
The famous Connes formula defines a spectral distance $d_s^\tau$ on
$\Xi$ from the spectral triple. The infimum
$\underline{d}_s:=\inf_\tau d_s^\tau$ coincides with $d$. We find the
following characterization of the tiling:
under certain assumptions on the sequence $(\delta_n)_n$
the supremum
$\overline{d}_s:=\sup_\tau d_s^\tau$ is Lipschitz-equivalent to
$\underline{d}_s$, i.e.\   $\exists c>0: \overline{d}_s\leq c
\underline{d}_s$, if and only if the continued fraction expansion of
$\theta$ is bounded. Thus our family of spectral triples is sensitive
to the arithmetic properties of the irrational $\theta$. It is known
that these properties are related to order: the continued fraction 
expansion of $\theta$ is bounded if and only if the tiling %we started out with 
is {\em linearly repetitive}, a property commonly regarded as the strongest notion of aperiodic order.
%which is commonly accepted as the strongest notion of aperiodic order.

%%%%%%%%%%%%%%%%%%%%%%%%%%%%%%%%%%%%%%%%%%
\subsection*{Summary of results}
We build spectral triples for compact metric spaces $(X,d)$, in particular for ultra metric spaces. 
Similar constructions have already been proposed in the literature:
see for instance \cite{Ri04, CI07, Pal10} for general metric spaces,
and \cite{La97,GI03,GI05,PB09,JS10a} for the case of fractals or
tiling spaces.
We believe that our construction is more natural.
In this paper we mainly concentrate our analysis on compact ultra metric spaces, and represent $X$ first, as in \cite{PB09}, by its Michon tree of clopen partitions \cite{Mich85}.
For a unilateral subshift the Michon tree is its {\em tree of words} \cite{Ca97} and for a tiling its {\em tree of patches}.
%We believe that our construction is more natural.
%%As in \cite{PB09} and similar to \cite{Pal10} the space $X$ is
%The space is represented by a tree $\Tt = (\Tt^{(0)}, \Tt^{(1)})$. 
%In the ultra metric case $\Tt$ corresponds to the Michon  tree of
%clopen partitions of $X$ \cite{Mich85}, which for a unilateral
%subshift is its {\em tree of words} \cite{Ca97} and for a tiling its
%{\em tree of patches}. 
%% a tree graphwhose set of infinite paths $\Pi_\infty$ is dense in $X$ or even homeomorphic to $X$ in some cases of interest.
We also use choice functions although in a slightly different way than in \cite{PB09}. 
Indeed, for us a choice function $\tau$ is used to define another graph, the approximating graph
$\Gamma(\tau)$, whose vertices form a dense subset of $X$ and whose edges encode which points are 
considered to be neighbors. We give each edge a length, namely the distance of the neighbors measured with the metric $d$.
This way $\Gamma(\tau)$ becomes a metric graph.  
%The graph metric $d_g$ on 
%$\Gamma(\tau)$ hence induces a metric on that dense set $V$.

Using the graph $\Gamma(\tau)$, we build a spectral triple for the $C^\ast$-algebra $C(X)$.
% with the sup norm, for each choice function $\tau$.
The Connes distance, or spectral distance $d_s^\tau$, defined by the spectral triple plays a crucial role in our characterization of aperiodic order. To obtain choice independent quantities we look at 
\( \underline{d}_s = \inf_\tau d_s^\tau \) and at  \( \overline{d}_s = \sup_\tau d_s^\tau \).
We have the following basic facts: 
%%%%%%%%%%%%%%
\begin{itemize}
\item $d_s^\tau$ is an extension of the graph metric.
%\item[$\circ$] $d_s^\tau$ is an extension of the graph metric.
%$d_g$ on $X$, for each $\tau$;
%\item $d \le d_s^\tau$, for each $\tau$;
%\item[$\circ$] \( \underline{d}_s  = d\).
\item \( \underline{d}_s  = d\).
\end{itemize}
%%%%%%%%%%%%%%
We ask the following questions:
%%%%%%%%%%%%%%
\begin{enumerate}[(a)]
\item When is $\overline {d}_s$ continuous and hence
induces the same topology?\label{item-a}
%\item $\bar{d}_s$ is continuous w.r.t. to $d$,
%\item When is $\bar d^\tau_s$ Lipschitz equivalent to $d$ (\(d \le d^\tau_s\le c d\), for some $c>0$),
\item When is $\overline{d}_s$ Lipschitz equivalent to $d$ (\(d \le  \overline{d}_s \le c d\), for some $c>0$)?\label{item-b}
\end{enumerate}
%%%%%%%%%%%%%%
The answer to the first question is not always positive. This is what we want. Thus we may characterize a metric space $(X,d)$ according to whether it satisfies (\ref{item-a}) or (\ref{item-b}) or neither.
We advocate that for discrete tiling spaces to satisfy (\ref{item-b}) is a {\em characterization of high order}. 
\medskip

We study in detail two examples related to aperiodic order. More specifically, we consider the discrete tiling space arising from
%%%%%%%%%%%%%%%%%%%
\begin{enumerate}[(1)]
%\begin{itemize}
\item \label{item-1}
a repetitive one-sided subshift which, for any $n\in \NM$, has a single right special word of length $n$,
\item \label{item-2}
a tiling of $\RM^d$ with finite local complexity.
%\end{itemize}
\end{enumerate}
%%%%%%%%%%%%%%%%%%%
Subshifts can be seen as tilings of $\R$ but we have more structure in this case and hence stronger results. 
% and  correspond to tilings of $\RM$ obtained by the canonical projection method with codimension $n$.
Discrete tiling spaces are known to be metrizable by an ultra metric. But there is no a priori best choice for this metric. We consider therefore a whole family of metrics on these tiling spaces and try to answer the above questions for all of them. This way we get a characterization of the tiling space, and hence of the tiling, 
but we see also a sign which singles out the most commonly used choice for the metric, namely that the distance between two tilings is the inverse of their coincidence radius.
%We provide necessary and sufficient conditions for this to happen in the case (\ref{item-1}), and a necessary condition in the case (\ref{item-2}).

The elements of a subshift are sequences of letters, which in the one-sided case are functions from $\N$ into some (finite) alphabet. 
The allowed sequences are not arbitrary but constrained and in our cases the constraints are such that all words (finite strings appearing in the sequence) of  given length $n$, except one, have a unique extension by a letter at their right end.
% (the exception is called the right special word of length $n$). 
Episturmian subshifts \cite{GJ09}, such as Arnoux-Rauzy subshifts for instance 
\cite{AR91}, are examples of such subshifts.
We say that the subshift is {\em repulsive} if there is a constant $\ell>0$ such that whenever a word of length $n$ occurs both as a prefix and as a suffix in a word of length $N>n$ then $\frac{N-n}{n}\geq \ell$.
A subshift is repulsive if and only if it is {\em power-free}, that is, it does not contain arbitrarily large powers of words.

Given a strictly decreasing sequence of positive real numbers $(\delta_n)_{n\in\NM}$ converging to zero, we define a ultra metric $d$ on the space of our sequences $\Xi$ by
%%%%%%%%%%%%%%%%%%
\[
d(\xi, \xi') = \inf \big\{ \delta_n \, : \, \xi_m = \xi'_m, \forall m \le n \big\}\,.
\]
%%%%%%%%%%%%%%%%%%
Our main result for subshifts is the following.

\vspace{.1cm}
\noindent {\bf Theorem \ref{thm-epistur}} 
{\it Consider a repetitive one-sided subshift
which has exactly one right special word per length. %as in (1). 
Let  $(\delta_n)_n$ be a (strictly decreasing) null-sequence such that 
\(\underline{c}\, \delta_n\leq\delta_{2n}\) and 
\(\delta_{nm}\leq \overline{c}\, \delta_n\delta_m\) 
for some $\overline{c},\underline{c}>0$. 
We provide the subshift space with the metric $d$ defined by $(\delta_n)_n$.
The following are equivalent:
\begin{enumerate}
\item   The subshift is repulsive.
\item   $\overline{d}_s$ is Lipschitz equivalent to $d$.
\end{enumerate}
}

%The {\em right-special repulsiveness} property is a measure of order of subshifts satisfying (2) (see Section~\ref{}).
A linearly repetitive (or linearly recurrent) subshift %\cite{DHS99,Du00} 
is repulsive, and both notions are equivalent for Sturmian subshifts (see Lemma~\ref{lem-repulsive}, and Section~\ref{tilings} for the definitions).
Thus the above specializes to the following.
%The {\em right-special repulsive} property is weaker than linear recurence, and is equivalent to it for Sturmian subshifts, hence the following.

\vspace{.1cm}
\noindent {\bf Corollary \ref{cor-stur}} {\it Consider a Sturmian subshift associated with the 
irrational \(\theta \in (0,1)\).
The following are equivalent:
%%%%%%%%%%%%%%
\begin{enumerate}[(i)]
\item The Sturmian sequence is linearly recurrent;
\item The continued fraction expansion of $\theta$ is bounded;
\item $\bar{d}_s$ is Lipschitz equivalent to $d$.
%\item $\theta$ has bounded continued fraction expansion.
\end{enumerate}
%%%%%%%%%%%%%%
}
Note that $\delta_n=e^{-n}$ does not satisfy the hypothesis of Theorem~\ref{thm-epistur}.
(The cases when the $\delta_n$'s decay exponentially fast are not interesting because then $\bar{d}_s$ and $d$ are always Lipschitz equivalent, see Theorem~\ref{thm-exp}). 
Possible choices for the $\delta_n$ satisfying the hypothesis are of the form
\(\delta_n = \frac{\ln^bn}{n^a}\) for $a>0$, $b\geq 0$.

On the other hand, continuity of $\overline{d}_s$ is always guaranteed if the series $\sum_n\delta_n$ is summable (see Corollary~\ref{cor-cont}).
We can show further that there exist Sturmian subshifts for which (a) and (b) fail when $\sum_n\delta_n$ is not summable: $\bar{d}_s$ is not even continuous (see Theorem~\ref{thm-cexstur}). 
In the case when $\delta_n= 1/n^a$, for $a>1$ the metric $\bar{d}_s$ is always continuous, and for $0<a\le 1$ there are Sturmian subshifts for which it is not.
This can be seen  as a sign that the metric corresponding to $a=1$, so given by $\delta_{n} = 1/n$, has something particular to it.
\medskip

We then consider higher dimensional tilings with finite local complexity. 
Here we have less structure
and so obtain only a sufficient condition for (\ref{item-b}).  We study its discrete tiling space $\Xi$ 
which is the set of all tilings (with a punctured tile at the origin) 
whose patches are copies of patches of $T$. We equip it
with a metric of the following type. Let $r(\xi,\xi')$ be the coincidence radius of two tilings $\xi,\xi'\in\Xi$, that is, the radius of the largest patch around the origin $0\in\R^d$ on which they coincide. Given any strictly decreasing function $\delta$ tending to zero at infinity, we define 
%%%%%%%
\[
d(\xi,\xi') = \delta(r(\xi,\xi')).
\]
%%%%%%%
The tiling is said to have {\em equidistributed frequencies} if the ratio of the smallest of the frequencies of patches of radius $r$, $\textrm{freq}_{\textrm{min}}(r)$, and the  largest one, $\textrm{freq}_{\textrm{max}}(r)$, satisfies:
%%%%%%%
\[
\frac{\textrm{freq}_{\textrm{min}}(r)}{\textrm{freq}_{\textrm{max}}(r)} \ge c \,, \quad \text{\rm for some constant $c>0$, and all $r>0$.} 
\]
%%%%%%%
%have the same asymptotic behavior as $r\rightarrow +\infty$.
For example, if the tiling is {\em linearly repetitive}, then it has equidistributed frequencies, and one has \( cr^{-d} \le \textrm{freq}_{\textrm{min}}(r) \le \textrm{freq}_{\textrm{max}}(r) \le C r^{-d}\) %where $\sim$ is understood up to a constant factor which may differ for the minimal and the maximal frequency 
(see Theorem~\ref{thm-freq}).
Our main theorem for tilings is the following.

\vspace{.1cm}
\noindent {\bf Theorem \ref{thm-Lip}} {\it Consider a tiling $T$ of $\RM^d$ with finite local complexity and equidistributed frequencies.
%Let $\Xi$ be its canonical transversal.
Suppose that the patch counting function $P$ satisfies $P(4r)\leq \tilde c P(r)$ for some $\tilde c$.
If the function $\delta$ satisfies 
%%%%%%%%%%%%%%
\begin{enumerate}[(i)]
\item \(\delta(ar) \le \delta(a) \delta(r)\), for all \( a, r \ge 1\), and  
\item \( \delta\in L^1 \big( [1,\infty), \frac{dx}{x} \big)\),
\end{enumerate}
%%%%%%%%%%%%%%
then $\overline{d}_s$ is Lipschitz equivalent to $d$.}
\medskip

There is a natural zeta-function associated to our spectral triples: it does not depend on the choices $\tau$, and is given by a series of powers of the eigenvalues of the Dirac operator.
It is similar to those obtained in \cite{PB09,JS10a}.
In particular for tilings which have a uniform bound on the number of possible patch extensions the abscissa of convergence of those various zeta-functions are all the same (see Lemma~\ref{lem-PB-zeta}).
In \cite{PB09} the abscissa of convergence %of the zeta-function 
was proven to be a fractal dimension of the space (namely, the upper box dimension).
In \cite{JS10a}, in the case of substitution tilings, it was identified with the exponent of complexity of the tiling, and further with the Hausdorff dimension of its discrete tiling space in \cite{JS10c}.
Our result here relates the abscissa of convergence in general to the weak complexity exponent. 
We define the lower and upper complexity exponents
%%%%%%%%%
\[
\underline{\beta} = \sup\{\gamma : P(r)\geq  r^\gamma, r\mbox{ large}\}\,, \quad
 \overline{\beta} = \inf\{\gamma : P(r)\leq r^\gamma, r \mbox{ large}\}\,,
\]
%%%%%%%%%
where $P$ is the patch counting function of the tiling or the subshift.
If \(\underline{\beta} = \overline{\beta}\) we call it the {\em weak complexity exponent}.
We present general theorems for tilings and subshifts, an example of which is for instance the following.

\vspace{.1cm}
\noindent {\bf Theorem \ref{thm-zeta-2}} (A particular case)
{\it For a $d$-dimensional tiling with finite local complexity and
\(\delta(r) \in L^{1+\epsilon}\big([1,+\infty)\big) \setminus L^{1+\epsilon}\big([1,+\infty)\big)\) for all $\epsilon>0$, the abscissa of convergence $s_0$ of the zeta function satisfies
%%%%%%%%%
\[
\underline{\beta}\le s_0 \le \overline{\beta} + d-1\,.
\]
%%%%%%%%% 
If the tiling has a uniform bound on the number of possible patch extensions and 
\(\underline{\beta} = \overline{\beta}\) then $s_0$ coincides with the weak complexity exponent.}
\medskip

We also derive Dirichlet forms and Laplace-Beltrami like operators from our spectral triples, and we compare those with the constructions given in \cite{PB09,JS10a} (see Section~\ref{sec-comparison}).
Our constructions  recover those previous results as particular cases.
There is a close relation with the Pearson-Bellissard Laplacians and even equality with
the zeta-function when the vertices of the Michon tree branch in at most two points.
But if this is not the case, the Pearson-Bellissard approach does not produce a metric when taking the supremum over all $d_s^\tau$.
In particular, one cannot use it for the characterization of order as we do.

%%%%%%%%%%%%%%%%%%%%%%%%%%%%%%%%%%%%%%%%%%%%%%%%%%%%%%%%%%%%%%%%%%%%%%%%%%%%
\subsection*{Aknowledgements}
This work was supported by the ANR grant {\em SubTile} no. NT09 564112.
%J.S. aknowledges postdoctoral funding from this grant.
The authors would like to thank L. Zamboni for helpful discussions and explaining them the relation between the various notions of repulsiveness used here (see Lemma~\ref{lem-repulsive}).

%%%%%%%%%%%%%%%%%%%%%%%%%%%%%%%%%%%%%%%%%%%%%%%%%%%%%%%%%%%%%%%%%%%%%%%%%%%%
\section{Notions from aperiodic order}

We present three basic notions of order. These can be applied to one-sided subshifts, $\Z^d$-subshifts or tilings of finite local complexity. We call them all simply tilings.

%%%%%%%%%%%%%%%%%%%%%%%%%%%%%%%%%%%%%%%%%%%%%%%%%%%%%%%%%%%%%%%%%%%%%%%%%%%%
\subsection{Subshifts and tilings}
\label{tilings}

Let $\A$ be a finite set of symbols, or letters.
On the space $\A^\N$ of one-sided infinite sequences (also called one-sided infinite words) we consider the left shift 
$\sigma$, defined by taking away the first letter.
Viewing $\A^\N$ as an infinite Cartesian product of copies of $\A$ we provide it with  the product topology. It is thus compact and the left shift is continuous.   Let
$\L$ be a collection of (finite) words in $\A$. The subshift space defined by $\L$, 
denoted $\Xi_\L$,  is the  subset of $\A^\N$ of all infinite sequences whose 
finite subwords belong to $\L$. 
We suppose that $\L$ contains all finite words of all infinite sequences of $\Xi_\L$. 
In fact, in our exemples below $\L$ is defined by a given sequence, namely as the set of all finite words of that sequence. The system $(\Xi_\L,\sigma)$ is called a {\em one-sided subshift}, and $\Ll$ its {\em language}.

We may as well consider two-sided infinite sequences $\A^\Z$ on which the left shift becomes bijective and hence, again in the product topology, a homeomorphism. $\Xi_\L$ is defined in the obvious analogous way and now we speak of a {\em two-sided subshift}. 
More generally, we may consider $\A^{\Z^d}$ with an action $\sigma$ of $\Z^d$ by coordinate wise left shift. One can picture these as decorations of the points of $\Z^d$ with symbols. Now $\L$ ought to be a collection of finite decorated subsets of $\Z^d$ and the space $(\Xi_\L,\sigma)$ is called a $\Z^d$-subshift.

A tile of $\R^d$ is a compact subset of $\R^d$ which is homeomorphic to a closed ball. 
A tiling of $\R^d$ is a covering by countable collection of tiles which do not overlap, except possibly at with their boundaries.  
One may also allow for decorations on the tiles and then an element of a $\Z^d$-subshift space can be seen as a tiling by decorated cubes. Such a tiling is called a Wang tiling.

A punctured tile is a tile together with a chosen point in its interior. 
A punctured tiling is a tiling whose tiles are punctured in a translationaly coherent way.
That is, if two tiles $t,t'$, with punctures $x,x'$, are translate of each other (possibly with matching collars): $t=t'+a$, for some $a \in \RM^d$, then so are their punctures: $x=x'+a$. 
If $T$ is a punctured tiling, we denote by $T^\punc$ the set of punctures of the tiles of $T$.

\medskip

We need to define the concept of an $r$-patch where $r$ is a positive real number.
For a one-sided subshift this is most naturally done using the length of words.
The length of a word $w$ is simply the number of letters of the word, which we denote by $|w|$.
Then an $r$-patch of the subshift is simply a word of length $[r]$ (integer part of $r$).

For a $\Z^d$ subshift an $r$-patch is a decorated cube of odd diameter
$l$ such that $\frac{l+1}{2} = [r]$. This seems a little artificial
but for comparison with tilings we want that these cubes are centered
around the center of a cube of edge-length $1$. 

For a punctured tiling $T$ of $\RM^d$, an {\em $r$-patch} is the set
$B_{0,r}[T-x]$ of tiles in $T-x$, for some $x\in T^\punc$, whose
punctures lie in the closed $r$-ball about the origin.  
For any $x$, the function $r\mapsto B_{0,r}[T-x]$ is piecewise
constant (when given the set of $r$-patches the discrete topology) and  
there exists a strictly increasing sequence \((r_n)_{n\in \NM}\) of positive
real values such that $\{r_n:n\in \N\}$ corresponds to the collection
of all discontinuities of these functions for different choices of
$x\in T^\punc$. We call \((r_n)_{n\in \NM}\) the sequence of sizes of $T$.

% The set of all $r$-patches of $T$ is denoted $\Pp_r(T)$.
% There exists a non decreasing sequence \((r_n)_{n\in \N}\) of positive
% reals, such that 
The set $\Pp(T)$ of all patches of $T$ is thus the union of the 
$\Pp_{n}(T)$, the set of all $r_n$-patches, over $n\in \NM$.

%The {\em continous tiling space} $\Omega$ of $T$, is the collection of all tilings $T'$ whoses patches are patches of $T$: \(\Pp(T')=\Pp(T)\).
The {\em discrete tiling space} $\Xi$ of $T$ is the set of tilings $T'$, such that \( \Pp(T')=\Pp(T) \) and \(T'^\punc \ni 0_{\RM^d} \) 
%%%%%%%%%%%%%%%%
\begin{equation}
\label{Xi}
\Xi = \big\{ T' \ \text{\rm tiling of } \RM^d \, : \,  \Pp(T')=\Pp(T)\,, \ T'^\punc \ni 0_{\RM^d} \big\}.
\end{equation}
%%%%%%%%%%%%%%%%
If one does not require \(T'^\punc \ni 0_{\RM^d}\), then one defines the continuous tiling space $\Omega$ of $T$.
For some classes of tilings $\Omega$ can alternatively be defined as the closure of the $\RM^d$-orbit of $T$ under translation in a certain topology for which it is compact, Hausdorff, and metrizable.
In this case $\Omega$ is an $\RM^d$-dynamical system, and $\Xi$ is a {\em canonical transversal} for the $\RM^d$-flow.

The following is well-known \cite{Kel95}.
%%%%%%%%%%%%%%%%
\begin{prop}\label{prop-metric}
The discrete tiling space is compact, Hausdorff, and metrizable. 
Any strictly decreasing function $\delta$ on $\R^+$ which tends to $0$ at infinity defines a compatible metric by 
%%%%%%%%%%%%%%%%
\[
d(\xi,\eta) = \inf \big\{ \delta(r) \, : \, B_{0,r}[\xi]=B_{0,r}[\eta]\} \,,
\]
%%%%%%%%%%%%%%%%
where $B_{0,r}[T]$ denotes the $r$-patch of $T \in \Xi$ centered at the origine.
Moreover, $d$ is a ultra metric, that is one has
%%%%%%%%%%%%%%%%
\begin{equation}
\label{ultrametric}
d(x,y) \leq \max\{d(x,z),d(y,z)\}\,, \quad \forall x,y,z \in \Xi.
\end{equation}
%%%%%%%%%%%%%%%%
\end{prop}
%%%%%%%%%%%%%%%%

A Delone set of $\RM^d$ is a discrete set of points for which the distances from any point to any of its nearest neighbors is uniformly bounded.
Namely, given $0< \rd < \ru$, a discrete set $\Ll$ of $\RM^d$ is a $(\rd,\ru)$-Delone set if it is both
%%%%%%%%%%%%%%%%
\begin{enumerate}[(i)]
%\begin{itemize}
\item $\rd$-{\em uniformly discrete}: any ball of radius $\rd$ intersects $\Ll$ in {\em at most} a point;

\item $\ru$-{\em relatively dense}: any ball of radius $\ru$ intersects $\Ll$ in {\em at least} a point;
%\end{itemize}
\end{enumerate}
%%%%%%%%%%%%%%%%
Given a Delone set, there is a natural associated punctured tiling, called its {\em Voronoi tiling}.
The tiles are convex polytopes: the {\em Voronoi cells} $v_x$ punctured at the points $x$ in $\Ll$:
%%%%%%%%%%%%%%%%
\begin{equation}
\label{voronoicell}
v_x = \big\{ y \in \RM^d \, : \, \|x-y\| \le \| z - y \|, \forall z \in \Ll \big\}\,.
\end{equation}
%%%%%%%%%%%%%%%%

%%%%%%%%%%%%%%%%%%%%%%%%%%%%%%%%%%%%%%%%%%%%%%%
\subsection{Patch counting and complexity}
%%%%%%%%%%%%%%%%
\begin{defn}
The patch counting function $P:\R^+\to \N$ of a tiling assigns to $r$ the number of $r$-patches modulo translation which occur in the tiling. $P$ is sometimes also called complexity function.
%In other words
%%%%%%%%%%%%%%%%
%\[
% P(r) = \# \big\{ B_r[\xi] \, : \, \xi\in\Xi \big\}\,,
%\]
%%%%%%%%%%%%%%%%
\end{defn}
%%%%%%%%%%%%%%%%
We suppose that this function is finite for all $r$, i.e.\ the tiling has {\em finite local complexity}. 
One expects that an ordered tiling has a small patch counting function, one that certainly should be
sub-exponential or even bounded by a polynomial function. 
We call
\begin{eqnarray*}
 \underline{\beta} &=& \sup\{\gamma: P(r)\geq r^\gamma \mbox{ for large $r$}\}\\
 \overline{\beta} &=& \inf\{\gamma: P(r)\leq  r^\gamma \mbox{ for large $r$}\}
 \end{eqnarray*}
the lower an upper complexity exponent of the tiling. 
Note that these quantities may be alternatively obtained as
 $\overline{\beta} = \overline{\lim}_{r\to \infty}\frac{\ln P(r)}{\ln r}$ and
 $\underline{\beta} = \underline{\lim}_{r\to \infty}\frac{\ln P(r)}{\ln r}$.
If both exponents are equal we call their common value the {\em weak complexity exponent} of the tiling.
Note that if  there exist constants $c,C>0$ such that
%%%%%%%%%%%%%%%%
\(
c r^{\beta} \leq P(r) \leq C r^{\beta} \,
\)
then $\beta$ is the weak complexity exponent. The latter behaviour for patch counting functions was studied in great detail in \cite{Ju09,Ju10} and there called polynomial growth.
Note that $\overline{\beta}=\underline{\beta}$ does not imply polynomial growth, as for instance $P(r) = r\ln r$ satisfies $\overline{\beta}=\underline{\beta}=1$.
Another stronger definition  of complexity exponent (which we do not investigate) would be to demand that $P(r)$ and $r^\beta$ are asymptotically equal, i.e.\ 
\({\lim}_{r\to \infty}\frac{P(r)}{r^\beta}=1\). 

%%%%%%%%%%%%%%%%
%If $\beta = 0$ then the tiling is periodic.
%It seems plausible to consider that $\beta = d$ is minimal for aperiodic repetitive tilings, but this is actually not the case \cite{Ca-private}.

The patch counting function is integer valued  and right continuous.
It has discontinuities at points which form a subsequence of the
sequence of sizes $(r_n)_{n\in\NM}$ introduced in Section~\ref{tilings}.

%%%%%%%%%%%%%%%%%%%%%%%%%%%%%%%%%%%%%%%%%%%%%%%%%%%%%%%%%%%%%%
\subsection{Repetivity}
%%%%%%%%%%%%%%%%
\begin{defn}
The repetivity function $R:\R^+\to\R^+$ of a tiling assigns to $r$ the smallest $r'$ such that any $r'$-patch of the tiling contains any $r$-patch up to translation.
\end{defn}
%%%%%%%%%%%%%%%%
We suppose that this function is finite for all $r$, i.e.\ the tiling is {\em repetitive}. Again, 
an ordered tiling is expected to have small repetitivity function. There are various bounds between the repetitity function and the complexity function. For example, Lagarias and Pleasants obtain 
%%%%%%%%
\[
R(r)  \geq c P(r)^\frac{1}{d} 
\]
%%%%%%%%
with some $c>0$ \cite{LP03} for repetitive $d$-dimensional tilings (and large $r$)\footnote{Lagarias and Pleasants work actually with Delone sets and the constant $C$ is related to their inner radius $\underline{\rho}$, but the results carry over to tilings.}.
A tiling is called {\em linearly repetitive} \cite{DHS99,Du00,LP03} if 
%%%%%%%%
\[
 R(r) \leq \clr \, r
\]
%%%%%%%%
for some $\clr>0$. 
In the context of subshifts, the terminology {\em linearly recurrent} is
traditionally used instead of linearly repetitive. 
Lenz shows in \cite{Le04} that for linear repetitive tilings there exists a constants $C>0$ and 
$\clr>0$ such that the reverse inequalities 
$$R(r)  \leq C P(r)^\frac{1}{d}  $$ and
$$ R(r) \geq \clr\,r $$
hold as well. The latter is actually a consequence of the repulsion property which we discuss below (see Theorem~\ref{thm-Lenz}). 
It follows that linearly repetitive tilings have complexity exponent $d$.
Linearly repetitive tilings are considered to be the most ordered
aperiodic tilings. 
\subsection{Repulsiveness for subshifts}
\begin{defn} A subshift is repulsive if
%%%%%%%%%%%%%%%%
\[
\ell = \inf\left\{\frac{|W|-|w|}{|w|} : w,W\in\L, \mbox{$w$ is a suffix and a prefix of $W$}, w\neq W\right\}
\]
%%%%%%%%%%%%%%%%
is strictly greater than $0$.

A word $w\in\L$ is called right special if $wa,wb\in \L$ for at least two distinct letters $a,b$.
We denote by $\L_{\Rr}$ the set of all right special words of $\L$.
We say that the subshift is right special repulsive if
%%%%%%%%%%%%%%%%
\[
\ell_\Rr = \inf\left\{\frac{|W|-|w|}{|w|} : w,W\in\L_\Rr, \mbox{$w$ is a suffix and a prefix of $W$}, w\neq W\right\}
\]
%%%%%%%%%%%%%%%%
is strictly greater than $0$.
\end{defn}

%%%%%%%%%%%%%%%%
\begin{lemma}
\label{lem-repulsive}
An aperiodic subshift is right special repulsive if and only if it is repulsive.
\end{lemma}
%%%%%%%%%%%%%%%%
{\em Proof:}
First note that \( \ell \le \ell_\Rr\) since the infimum is taken over a larger set for $\ell$, and therefore \( \ell_\Rr= 0 \Rightarrow \ell =0\).

To prove the converse, assume that $\ell=0$, and fix $0<\epsilon < 1/2$.
By definition of $\ell$ there are distinct words $w,W$, such that $w$ is a prefix and suffix of $W$, and one has \( (|W|-|w|) / |w| \le \epsilon/4\).
If $w$ is right-special, then $\ell_{|w|} \le \epsilon$, and therefore $\ell_\Rr \le \epsilon$.
% and since $\epsilon$ was arbitrary this proves $\ell_\Rr=0$.
If $w$ is not right-special, then we claim that there is an integer $j\le |W|-|w|$ such that $W_j=W_{[0,|w|+j]}$ (and therefore $w_j=w_{[0,|w|-j]}$) is right-special.
For if this were not true, then the right extension of $w$ up to $W$ would be forced, and since $w$ is also a suffix of $W$, thus also the right extension of $W$ would be forced, and we would inductively build an infinite periodic sequence this way.

Now assume that  $W_j$ and $w_j$ are right-special for some $j\le |W|-|w|$.
We have \( |W_j| - |w_j| = |w| + j - (|w|-j) = 2 j \le 2 (|W| - |w|)\),
and also \(|w_j| = |w|-j \ge (1-\epsilon) |w| \ge |w|/2\).
Therefore \( (|W_j| - |w_j|) /|w_j| \le 4   (|W| - |w|) /|w| \le \epsilon\).
Since $w_j$ is right-special this implies that \(\ell_\Rr \le  \epsilon\), and as $\epsilon$ was arbitrary this proves $\ell_\Rr = 0$.
\qed
\medskip

It may be shown that repulsiveness is equivalent to {\em power free}, i.e.\ the language 
does not contain arbitrary large powers of words.
%Equivalently, a subshift is {\em not} power free, if for all $n\in \NM$ there exists a word $u$ for which $u^n$ belongs to its language.

%%%%%%%%%%%%%%%%%%%%%%%%%%%%%%%%%%%%%%%%%%%%%%
\subsection{Tilings with equidistributed frequencies}
If all $r$-patches had equal frequency then this frequency ought to be
$P(r)^{-1}$. Hence we define 
%%%%%%%%%%
\begin{defn}
A tiling has {\em equidistributed frequencies} if there
are constants $c,C>0$ such that for all $r$-patches $p$
%%%
\[
 c P(r)^{-1} \leq \mbox{freq}(p) \leq C P(r)^{-1} .
\]
%%%
\end{defn}
%%%%%%%%%%
Let $p$ be an $r$-patch of a tiling $T$. 
We denote by  $\L_p$ the Delone set of occurrences of $p$ in $T$,
%%%
\[
\Ll_p = \big\{ x \in T^\punc \, : \, B_r[T-x] = p \big\}
\]
%%%
and let $\underline{\delta}(\L_p)$ and $\overline{\delta}(\L_p)$ be its parameters of uniform discreteness and relative denseness respectively.
For $y\in \L_p$ denote by $v_y$ its Voronoi cell in $\L_p$, and by $\ro(v_y),\ri(v_y)$ its outer and inner radii. 
$T$ satisfies the property of uniformity of return words if there exists $\sigma>0$ such that for all $r$-patches $p$ and $y \in \L_p$
%%%%%%%
\begin{equation}
\label{eq-U}
\ro(v_y) \le \sigma \ri(v_y).
\end{equation}
%%%%%%%
The reasonable generalization of repulsiveness for a tiling $T$ is that
\[
\ell := \inf\left\{\frac{\underline{\delta}(\L_p)}{r}:r>0,\mbox{$p$ $r$-patch of $T$} \right\}
\]
is strictly greater than $0$.

Repulsiveness thus says that the length of a return vector of a patch has to increase at least linearly with the size of the patch. Non-repulsiveness, i.e.\ $\ell(T)=0$ would imply that large patches overlap with themselves and larger and larger overlap occurs for larger and larger patches.

%\end{defn}
%%%%%%%%%%%%%%%%
\begin{thm} [\cite{Le04,BL08}]
\label{thm-Lenz} A linearly repetitive tiling is repulsive and satisfies (\ref{eq-U}).
Furthermore its patch frequencies exist and are given by  
%%%
\[
 \mbox{\rm freq(p)} = \lim_{n\rightarrow \infty} \frac{\#_p (C_n)}{vol(C_n)}\,.
\]
%%%
Here \(\#_p(C_n) = \#( \L_p \cap C_n)\) is the number of occurrences of $p$ in the cube $C_n$ of side-length $n$.
\end{thm}
%%%%%%%%%%%%%%%%
%%%%%%%%%%%%%%%%
\begin{thm}
\label{thm-freq}
A linearly repetitive tiling of $\RM^d$ has equidistributed frequencies: there exist constant  $c,C>0$, such that for all $r$-patches $p$ one has: 
%%%
\[
 c r^{-d} \leq \mbox{\rm freq}(p) \leq C r^{-d} \,.
\]
%%%
%Here $d$ is the dimension of the tiling. 
\end{thm}
%%%%%%%%%%%%%%%%
{\em Proof:}
%Let $\L$ be the set of punctures of the tilings which we suppose to be linearly repetitive.
Let $\omega$ denote the volume of the unit ball in $\R^d$. Let $p$ be an $r$-patch.
One has
%%%%%%%%%%%%%%%%%%
\[
\#_p (C_n) = \sum_{y\in C_n\cap \L_p} 1
\ge  \frac{1}{\omega} \sum_{y\in C_n\cap \L_p} \frac{vol(v_y\cap C_n)}{\ro(v_y)^d}
\ge \frac{1}{\omega \sigma^d} \sum_{y\in C_n\cap \L_p} \frac{vol(v_y\cap C_n)}{ \ri(v_y)^d}\,,
\]
%%%%%%%%%%%%%%%%%%
where the last inequality follows from (\ref{eq-U}).
Now we clearly have \(\ri(v_y) \le \overline{\delta}(\L_p)\), and, by linear repetivity  there exists a $C>0$ such that  \(\overline{\delta}(\L_p) \le C r\) ($C$ does not depend on $r$).
Hence \(\ri(v_y) \le C r \) and we get the lower bound
%%%%%%%%%%%%%%%%%%
\[
\#_p (C_n) \ge \frac{1}{\omega \sigma^d C^d} \sum_{y\in C_n\cap \L_p} \frac{vol(v_y\cap C_n)}{ r^d}
=  \frac{1}{\omega \sigma^d C^d} \, r^{-d}\; vol(C_n) \,.
\]
%%%%%%%%%%%%%%%%%%
Divide by $vol(C_n)$ and take the limit as $n\rightarrow \infty$.
The left-hand side yields $\freq(p)$ by Theorem~\ref{thm-Lenz}, and we get the lower bound of the statement.

For the upper bound one has
%%%%%%%%%%%%%%%%%%
\[
\#_p (C_n) = \sum_{y\in C_n\cap \L_p} 1
\le  \frac{1}{\omega} \sum_{y\in C_n\cap \L_p} \frac{vol(v_y)}{\ri(v_y)^d}
\le \frac{\sigma^d}{\omega} \sum_{y\in C_n\cap \L_p} \frac{vol(v_y)}{ \ro(v_y)^d}\,,
\]
%%%%%%%%%%%%%%%%%%
where the last inequality follows again from (\ref{eq-U}).
Now  \(\ro(v_y) \ge \underline{\delta}(\L_p)\), and by repulsiveness, \(\underline{\delta}(\L_p) \ge \ell r\).
Hence \(\ro(v_y) \ge \ell r\), and one has the upper bound
%%%%%%%%%%%%%%%%%%
\[
\#_p (C_n) 
\le \frac{\sigma^d}{\omega \ell^d} \sum_{y\in C_n\cap \L_p} \frac{vol(v_y)}{\underline{\delta}(P)^d} 
\le \frac{\sigma^d}{\omega \ell^d} \,r^{-d}\; \bigl( vol(C_n) + \overline{\delta}(\L_p)
area(\partial C_n) \bigr) \,,
\]
%%%%%%%%%%%%%%%%%%
where $area(\partial C_n)$ denotes the area of the boundary of the cube.
Now divide by $vol(C_n)$ and take the limit as $n\rightarrow \infty$.
The left-hand side yields $\freq(p)$ by Theorem~\ref{thm-Lenz}, while the term \(area(\partial C_n) / vol(C_n)\) goes to zero as $n$ tends to infinity and thus one gets
the upper bound of the statement.
\qed
\medskip

The converse of the theorem seems not to be true \cite{Lenz-private}.

%%%%%%%%%%%%%%%%%%%%%%%%%%%%%%%%%%%%%%%%%%%
\subsection{Uniform bounds on the number of possible patch extensions}
It seems that the following property of a tiling is of importance. 
Consider a tiling of finite local complexity. Then for each $r$-patch there are only finitely many possibilities of extending it to its neighborhood, one also says that there are only finitely many coronae for each $r$-patch. If the number of possibilities is uniformly bounded on the set of all $r$-patches we say that the tilings admits a uniform bound on the number of possible patch extensions.
  
A one-dimensional tiling and hence also a subshift satisfy this property trivially.
We will see later that tilings with equidistributed frequencies share this property if their patch counting function satisfy $P(4r) \le \tilde{c} P(r)$ (see Lemma~\ref{lem-equidist-UBPE}).
It would be interesting to have alternative characterisations.

%%%%%%%%%%%%%%%%%%%%%%%%%%%%%%%%%%%%%%%%%%%%%%%%%%%%%%%%%%%%%%%%%%%%%%%%%%%
\section{Spectral triples for compact metric spaces}

%%%%%%%%%%%%%%%%%%%%%%%%%%%%%%%%%%%%%%
\subsection{Background to spectral triples}
By a spectral triple for a compact (Hausdorff) space $X$ we mean a 
spectral triple for the algebra $C(X)$ of continuous functions over $X$.
It is hence given by a Hilbert space $\H$ with a faithful representation $\pi$
of $C(X)$ and a self-adjoint operator $D$, called the {\em Dirac operator}, satisfying: 
\begin{itemize}
\item[1.] the set $\tilde C(X)$ of continuous functions such that 
$[D,\pi(f)]$ extends to a bounded operator is a dense sub algebra of $C(X)$,
\item[2.] $D$ has compact resolvent.
\end{itemize}
Given such a spectral triple $(C(X),D,\H)$ the formula
\begin{equation}\label{eq-metric} d_s(x,y) :=
 \sup\{|f(x)-f(y)|: f\in C(X),\|[D,\phi(f)]\|\leq 1\} \end{equation} 
defines a pseudo-metric on $X$, called the spectral or Connes-metric. 
$d_s(x,y)$ may be infinite but if 
\begin{itemize}
\item[3.] the representation $\pi$ is non-degenerate and
\item[4.] the kernel of the derivation $[D,\pi(\cdot)]$ contains only constant functions
\end{itemize}
then $d_s$ is a metric \cite{Lat07}. It is, of course, important to
study the question when $d_s$ induces the topology on $X$. This turns out
to be equivalent to requiring that $\tilde C(X)/\C$ is pre-compact in the
quotient topology ($\C$ being the subspace of constant functions)
\cite{Rieffel}.  
While we will consider here triples which satisfy conditions $1.-3.$
and mostly also $4.$ we do not a priori require $d_s$ to induce the metric
but use this as a criterion for aperiodic order.
\bigskip

Given a pseudo-metric $d$ on $X$ the formula
\begin{equation}\label{equ-Lipschitz}
 \|f\|_{Lip,d} := \sup_{x\neq y}\frac{|f(x)-f(y)|}{d(x,y)}
\end{equation}
defines a seminorm on $C(X)$, called the Lipschitz-seminorm associated with $d$
(strictly speaking, $ \|f\|_{Lip,d}$ is finite only on a subset of $C(X)$ and so we allow here for $+\infty$
as a value). 
Given two pseudo metrics $d_1$ and $d_2$, it is easily verified that
\begin{equation}\label{eq-lip1}
\forall f : \|f\|_{Lip,d_1} \geq \|f\|_{Lip,d_2} \quad\mbox{iff}\quad \forall x,y:d_1(x,y)\geq d_2(x,y)  
\end{equation} 
so that the Lipschitz-seminorm determines uniquely its pseudo metric.

$\| [D,\pi(\cdot)] \|$ is also an example of a Lipschitz-seminorm\footnote{a seminorm which vanishes on constant functions and may take $+\infty$ as value.}
on $C(X)$ and a natural question is how it compares with $\|\cdot\|_{Lip,d_s}$. Rieffel shows that both are equal \cite{Rieffel}. While $\|f\|_{Lip,d_s}\leq \|[D,\pi(f)]\|$ follows easily
($\|[D,\pi(f)]\|\leq M$ implies $\forall x\neq y$: $d_s(x,y)\geq
M^{-1}|f(x)-f(y)|$ and the latter is equivalent to
$\|f\|_{Lip,d_s}\leq M$) the opposite inequality requires in general a lot more work.
In our context below it is also direct due to the particular form of the triple, see Lemma~\ref{lem-comparison}.

%%%%%%%%%%%%%%%%%%%%%%%%%%%%%%%%%%%%%
\subsection{Rieffel's triple}
Let $(X,d)$ be a compact metric space. Rieffel associated to this
space a triple, satisfying the axioms of a spectral triple except
possibly that of compactness of the resolvent.  
It is designed so that  its spectral distance reproduces $d$.
A modification of this triple will enforce compactness, at the price
of no longer reproducing exactly the metric.
This has been realized by a number of people now: see for instance \cite{Ri04, CI07, Pal10} for general metric spaces, and \cite{La97,GI03,GI05,PB09,JS10a} for the case of fractals. 
We describe here our own way of
looking at it. 

Let $\diag(X)$ be the diagonal in $X\times X$
and $\nu$ be a positive measure on
$Z=X\times X \backslash\diag(X)$ which is symmetric under the flip of
the points. Let 
$$\H = L^2(Z,\nu)$$
on which $C(X)$ acts via the representation $\pi$
\begin{equation}\label{eq-rep}
\pi(f)\psi(x,y) = f(x)\psi(x,y)
\end{equation}
and 
\begin{equation}\label{eq-Dirac}
D\psi(x,y) = d(x,y)^{-1}\psi(y,x).
\end{equation}
A direct calculation shows that 
$$\|[D,\pi(f)]\| = \sup\{\frac{|f(y)-f(x)|}{d(x,y)}:(x,y)\in \mbox{\rm
 supp}\nu\}.$$
If the support of $\nu$ is all of $Z$ then 
$\|[D,\pi(f)]\| = \|f\|_{Lip}$ and the usual calculus shows that  
the spectral distance coincides with $d$. The problem is that $D$ has
then rarely compact resolvent. 
So the trick is to adapt $\nu$ so that $D^{-1}$ becomes compact.

We start by decomposing $\H = \H_{+}\oplus \H_{-}$ according the
eigenspaces of the symmetry $\sigma(x,y) = (y,x)$. This symmetry
commutes with $D$ and so $D=D_+\oplus D_-$. In what follows we may
concentrate on $D_+$, as the analysis of $D_-$ is similar. 
It is then clear that $D_+^{-1}$ acts on $\H_+$ as left multiplication
by $d$ and hence has spectrum 
equal to the image of $d$. Thus for $D^{-1}$ to be compact we need to require:
\begin{enumerate}
\item the image of $d$ on $\supp\nu$ is discrete with perhaps one accumulation point at $0$,
\item  the subspace $L^2(Z_r,\nu)$ where
$Z_r = \{(x,y) : d(x,y)=r\}$  has to be finite dimensional for any $r>0$.
\end{enumerate}  
It follows that, if $X$ is not a finite set, then $\nu$ cannot have full support!
More precisely, for any compact $K\subset X\times X\backslash\diag(X)$ the intersection $\supp(\nu)\cap K$ must be a finite set.
In particular, $\nu$ is a countable pure point measure, i.e.\
%%%
\[
 \nu = \sum_{i\in I} {\nu_i}(\delta_{(x_i,y_i)}+\delta_{(y_i,x_i)})
\]
%%%
for some countable family of pairs in $X\times X$ indexed by $I$, $\nu_i\geq 0$.
The coefficients $\nu_i$ won't play a role in what follows and we will
set them equal to $1$. 

With these assumptions made we have a genuine spectral triple, i.e.\ a
Dirac operator with compact resolvent. Moreover, since $C(X)$ is
unital, the representation $\pi$ is non-degenerate.
We call this spectral triple the spectral triple associated to the measure $\nu$.

\begin{lemma}\label{lem-1}
In the above set-up  the spectral distance $d_s$ is larger than $d$.
\end{lemma}
{\em Proof:}
%The direct calculation above now shows that 
$$\|f\|_{Lip,d_s}\leq \|[D,\pi(f)]\| = \sup\{\frac{|f(y)-f(x)|}{d(x,y)}:(x,y)\in \supp(\nu)\}
\leq \|f\|_{Lip,d}.$$ 
%and since it is well-known that $\sup\{{|f(y)-f(x)|} : \|f\|_{Lip,d}\leq 1\} = d(x,y)$ the statement follows.
\qed
\begin{cor}
The spectral distance $d_s$ is compatible with the topology of $X$ if and only if it is continuous.
\end{cor}

\noindent{\em Proof:} 
Since $d\leq d_s$ 
%we have $d_s(x,y)<r\:\Longrightarrow d(x,y)<r$ and hence
%$B_{d_s}(x,r)\subset B_{d}(x,r)$. Hence 
any $d_s$-converging sequence converges also in the topology of $X$. Conversely let $(x_n)_n$ be a sequence converging to $x$. If $d_s$ is continuous then
$d_s(x_n,x)\to 0$. Hence $d_s$ is compatible with the topology on $X$. The converse is clear.  
\qed

%%%%%%%%%%%%%%%%%%%%%%%%%%%%%%%%%%%%%
\subsection{Christensen \& Ivan's sums of spectral triples for pairs}
The above described spectral triple can be understood as a countable
sum of spectral triples for pairs. 
This is the way Christensen {\it etal.}\ looked at it \cite{CI07}. One of
their question was, whether one can find a countable symmetric subset
of $X\times X$ such that
the above construction yields a spectral triple whose spectral distance reproduces the distance we started out with or, is at least as close as possible. Christensen {\it etal.}\ came up with existence results yielding a spectral distance which is Lipschitz equivalent with Lipschitz constant arbitrarily close to $1$, but not the same metric.
We sketch their construction in our own language.

\begin{enumerate}
\item Consider a sequence of $(\underline{\rho}_n, \overline{\rho}_n)$-Delone sets $\L_n\subset X$ such that the sequences $(\underline{\rho}_n)_{n}$ and $(\overline{\rho}_n)_{n}$ decrease to zero. (Their union $\bigcup_n \L_n$ could be understood as the
vertices of a graph.) 
\item Define horizontal pairs, that is, for each $n$ consider a symmetric subset
of points $E_n\subset \L_n\times \L_n\backslash \diag(\L_n)$. 
(These could be understood as unoriented edges of a graph.)
\item Define vertical pairs, i.e.\ for each $n$ a subset
 $I_n\subset\L_{n-1}\times\L_{n}$.
%These could be understood as oriented edges of a graph.
\end{enumerate}
This defines the measure
$$ \nu = \sum_n \left(\sum_{(x,y)\in E_n}  \delta_{x,y} +
 \sum_{(x,y)\in I_n}(\delta_{x,y}+\delta_{y,x})\right)$$ 
and hence a spectral triple by Rieffel's construction.
%%%%%%%
\begin{thm}[\cite{CI07}]
\label{thm-Christensen}
For any constant $C>1$ there exist a sequence of Delone sets
$\L_n\subset X$ together with a choice of horizontal pairs $E_n$ and
vertical pairs $I_n$ as above, such that the spectral triple
associated to the measure $\nu$ above yields a spectral metric $d_s$
which satisfies $$ d_s(x,y) \leq C d(x,y).$$
\end{thm}
%%%%%%%
For the proof we refer to the original article mentioning only that the 
%construction of the data producing the support of the measure is based on a sequence of $(\underline{\delta}_n,\overline{\delta_n})$-
Delone sets $\L_n$ are such that $\overline{\rho}_n$  decreases exponentially fast to $0$.
Moreover the distance between points in $E_n$ and $I_n$ is of the same order than $\overline{\rho}_n$.
It will become clear below that we are interested in a much slower decrease of $\overline{\rho}_n$ (see for instance Theorem~\ref{thm-exp}).
As a consequence, we will not be able to allow for Lipschitz constants $C$ which are arbitrarily close to $1$,
but this is not our issue.

%%%%%%%%%%%%%%%%%%%%%%%%%%%%%%%%%%%%%
\subsection{Our version} 
We slightly modify the above construction. In fact we demand that the
sequence of Delone sets is  
a chain, i.e.\ that $\L_n\subset \L_{n+1}$.  This way we do not need
vertical pairs, 
or one could also say that we chose vertical pairs of distance zero. Our
spectral triple is thus the given by the measure 
$$ \nu = \sum_n \sum_{(x,y)\in E_n}  \delta_{x,y} .$$
We still demand that $\bigcup_n\L_n$ is dense. 
We write $E$ for the union of the $E_n$.
By compactness, each $\L_n$ is finite. 
%Given the above graph approximation of $X$ we define
%$$\nu = \sum_n \sum_{e\in \E} \delta_e.$$
Rieffel's construction therefore provides us with the spectral
triple $(C(X),D,\ell^2(E))$, which for $e=(x,y) \in E$ reads
\begin{eqnarray}\label{eq-repr}
\pi(f)\psi(e) & = & f( s(e) ) \psi(e) \\
 D\psi(e) &=& l(e)^{-1}\psi(\tilde{e}). \label{eq-dirac} 
%\pi(f)\psi\big( (x,y) \big) & = & f(x) \psi\big( (x,y) \big) \\
 %D\psi\big( (x,y) \big) &=& l(x,y)^{-1}\psi\big( (y,x) \big). \label{eq-dirac} 
\end{eqnarray}
where $s(e)=x$ is the source of the edge $e$, and we have set $l(e) = d(x,y)$, and $\tilde{e}=(y,x)$ for the opposite edge to $e$.

A spectral triple is called an even spectral triple if there exists a $\Z_2$-grading of the Hilbert space, such that the elements represented by the algebra preserve the degree while  
$D$ exchanges the degree of vectors in the Hilbert space. This is the case with our spectral triples, though not in a canonical way.

Indeed, given that the edges appear in pairs, namely in two opposite  orientations, we devide them into two disjoint subsets, $E = E^+\cup E^-$, depending on a choice of orientation for each pair.
Having fixed such a choice the Hilbert space thus can be written 
$\ell^2(E) = \ell^2(E^+)\oplus \ell^2(E^-)$.
We see from \eqref{eq-repr} and \eqref{eq-dirac}) that $\pi(f)$ is indeed even, while $D$ is odd w.r.t.\ this grading.  

The zeta-function of our triple has the simple form
$$ \zeta(s) := \Tr(|D|^{-s}) = \sum_n\sum_{e\in E_n} l(e)^s.$$
%The factor $2$ arrises because $|D| = |D_+|\oplus |D_-|$ and $|D_-|=|D_+|$.

It is not so clear how to define a Dirichlet form on (real-valued)
continuous functions over $X$ 
using the commutator with the Dirac operator, that is, a quadratic
form of the type \((f,g)\mapsto \Tr(\rho(D)[D,\pi(f)][D,\pi(g)])\), where $\rho$ is a positive function (which we see as a ``density matrix'').
The difficulty lies in the correct choice for the Hilbert space. 
At first sight one would think that this ought to be the (real version of the) Hilbert space of the triple, but $f,g$ are continuous functions on $X$ and it seems difficult to relate them to $L^2$-functions on $Z$. 
We will see below that we can circumvent this difficulty by some sort
of average over the measures  
$\nu$.

%%%%%%%%%%%%%%%%%%%%%%%%%%%%%%%%%%%
\subsection{Metric graphs}
It should be clear by now that the above constructions are based on
an underlying graph. The vertices of the graph provide an
approximation of the space and the edges tell us which points  
are nearest neighbors. On the other hand, it is possible to define spectral
triples on metric graphs, {\it i.e.}\ unoriented graphs whose edges are
equipped with lengths.
It will be useful to develop this point of view. 

Let $\Gamma=(V,E)$ be a graph whose edges are provided
with lengths. $V$ is the set of vertices and $E\subset V\times V\backslash\diag(V)$, the
set of edges, a symmetric subset equipped with a function 
$l:E\to\R^+$ satisfying $l(x,y) = l(y,x)$.
{Our graphs thus fulfill the requirement that their are no loop edges and at most one edge between two vertices. Furthermore edges always occur in both
 directions, i.e.\ the graph is unoriented. 
We suppose also that $V$ is countable and the graph connected. Note that $V$ is in general not compact.

The length function $l$ defines a distance function on $V$, called the
graph metric $d_g$, namely $d_g(x,y)$ is the sum of the lengths of the
edges of the shortest path between $x$ and $y$:
%%%
\[
 d_g(x,y) = \inf \sum_{k=0}^N l((x_k,x_{k+1}) 
\]
%%%
the infimum running over all (finite) sequences $(x_k)_{0\leq k\leq N}$ such that
$x=x_0$, $x_N=y$, $(x_k,x_{k+1})\in E$.

% We are particularily interested in the case that the graph is bounded,
% that is, $d_g$ is finite, or even pre-compact, that is, for any
% $\epsilon$ there exists a finite $\epsilon$-cover. 

Rieffel's construction suggest to consider
the following triple over $C_u(V)$ the
algebra of uniformly continuous functions on $V$: 
The Hilbert space is $\H = \ell^2(E) = L^2(Z,\nu)$ where $Z=V\times
V\backslash\diag(V)$ and
$ \nu = \sum_{e\in E} \delta_e$, the representation as in
(\ref{eq-rep})
and the Dirac operator as in (\ref{eq-Dirac}). 
If for any given $r>0$ there are only finitely many edges of length greater than $r$ then $D^{-1}$ is compact.
The zeta-function is 
$ \zeta(s) = \sum_{e\in E} l(e)^s$
and the spectral distance on $V$ is the graph distance.

Note that the graph Laplacian (the metric one, not the combinatorial
one) can be obtained from  
the Dirichlet form 
$$ (f,g) \mapsto \Tr ([D,\pi(f)][D,\pi(g)]) $$
defined on $\ell^2(V,\R)$. But again, $C_u(V,\R)$ is rather different from
$\ell^2(V,\R)$.

%%%%%%%%%%%%%%%%%%%%%%%%%%%%%%%%%%%%%%%%%%%%%%%%%%%%%
\subsection{Densely embedded graphs}
\label{embdgraph}
%%%%%%%%
\begin{defn}
A dense embedding of a metric graph into a compact metric space is a
uniformly continuous injection $\tau:(V,d_g)\to (X,d)$ such that
$\tau(V)$ is dense.
\end{defn}
%%%%%%%%
Since $\tau^*(C(X))\subset C_u(V)$ we obtain a spectral triple $(C(X),\ell^2(E),D)$ over $C(X)$, with representation $\pi_\tau = \pi \circ \tau^\ast$ that is
%%%
\[
\pi_\tau(f)\psi(v,w) = f(\tau(v)) \; \psi(v,w) \,, %\quad D \psi(v,w) = l(v,w)^{-1} \psi(w,v)\,,
\]
%%%
and with $D$ given as in \eqref{eq-dirac} (or \eqref{eq-Dirac}).
This representation is faithful, since $\tau(V)$ is dense in $X$.
We are particularly interested in the case that the closure of 
$\tau\times\tau(E)$ in $X\times X$ contains the diagonal, in which
case $\tau(V)$ is dense in $X$.

We consider now the case in which  
%%%
\[
 d(v,w) = l(\tau(v),\tau(w))\,, \quad \forall (v,w) \in E\,.
\]
%%%
Uniform continuity of $\tau$ is then automatic, as
$d\left|_{V\times V}\right. \leq d_g$ by the triangle inequality. 
We wish to compare the metrics $d$, $d_g$, and the spectral distance
of the above triple which is given by
%%%
\[
d_s(x,y)  = \sup \{|f(x)-f(y)|:  {\forall (v,w)\in E}\;
|f(\tau(v))-f(\tau(w))|\leq l(v,w)\}.
\]
%%%
%%%%%%%%%%%
\begin{lemma}\label{lem-comparison}
$d_s$ is an extension of $d_g$, i.e.\ $ d_s\circ \tau\times\tau = d_g$. Furthermore.
$\|\tau^*f\|_{Lip,d_g}=\|[D,\pi_\tau(f)]\|=\|f\|_{Lip,d_s}$.
\end{lemma}
%%%%%%%%%%%
{\em Proof:} We first show the remaining inequality $\|f\|_{Lip,d_s}\geq \|[D,\pi_\tau(f)]\|$.
Let $(v,w)\in E$. Then
$d_s(\tau(v),\tau(w))=\sup_{f\in C(X)}\{|f(\tau(v))-f(\tau(w))|:
\forall (v',w')\in E: |f(\tau(v'))-f(\tau(w'))|\leq l(v',w')\} \leq l(v,w)$. 
%Hence by  $d\circ\tau\times\tau \leq d_g$ by the triangle inequality. 
Since also $d \leq d_s$ we have
$d_s(\tau(v),\tau(w))= l(v,w)$ for $(v,w)\in E$. Now
$\forall x\neq y$: $d_s(x,y)\geq
M^{-1}|f(x)-f(y)|$ implies 
%$\forall x\neq y\in V$: $d_g(x,y)\geq M^{-1}|f(x)-f(y)|$ implies 
$\forall (v,w)\in E$: $l(v,w)=d_s(\tau(v),\tau(w))\geq
M^{-1}|f(\tau(v))-f(\tau(w))|$ implies $\|[D,\pi_\tau(f)]\|\leq M$.
%We thus have shown that $\|f\|_{Lip,d_s}= \|[D,\pi(f)]\|=\|\tau^*f\|_{E-Lip} $.

To prove $\|[D,\pi(\tau^*f)]\| = \|\tau^*f\|_{Lip,d_g}$ note first that, if $x_0 < x_1 <\cdots x_f \in \R$
then, for any continuous function $f:\R\to\C$, $\sup_{k\neq j}\frac{|f(x_k)-f(x_j)|}{|x_k-x_j|} = \max_k
\frac{|f(x_k)-f(x_{k-1})|}{|x_k-x_{k-1}|}$. Hence if we consider two vertices $v,w$ in the graph and 
a shortest path $x_0,x_1,\cdots, x_f$ between them, then embedding the path isometrically into $\R$ shows that $\frac{|f(v)-f(w)|}{d_g(v,w)} \leq \max_k
\frac{|f(x_k)-f(x_{k-1})|}{l(x_{k-1},x_k)}$. This shows that $\|\tau^*f\|_{Lip,d_g}\leq \|[D,\pi(\tau^*f)]\|$.
The opposite inclusion is direct. Hence $\|\tau^*f\|_{Lip,d_g}= \|[D,\pi(\tau^*f)]\|$ which implies $ d_s\circ \tau\times\tau = d_g$.  
% In fact, $d(\tau(v),\cdot)$ is a continuous
% function which satisfies $|d(\tau(v),x)-d(\tau(v),y)|\leq d(x,y) \leq
% d_g(\tau^{-1}(x),\tau^{-1}(y))$  
% provided $v,x,y\in \tau(V)$. Hence for all $v$
%  \begin{eqnarray*}d_s(x,y)  &=& \sup \{|f(x)-f(y)|:  {\forall v,w,}\;
% |f(\tau(v))-f(\tau(w))|\leq d_g(v,w)\}\\
% & \geq& |d(v,x)-d(v,y)|
% \end{eqnarray*}
% By continuity of $d$ we have $\inf_v |d(v,x)-d(v,y)| = d(x,y)$.
\qed
\medskip

As a consequence of the last lemma, many properties of $d_s$ can be traced back to properties of $d_g$.
%%%%%%%%%%%
\begin{cor}
The spectral distance $d_s$ is bounded if and only if the graph is bounded.
A sufficient condition for this is that $d_s$ is continuous. 
In particular, if the graph is unbounded then $d_s$ cannot be continuous.  
\end{cor}
%%%%%%%%%%%
\noindent{\em Proof:} 
Suppose that $d_s$ is continuous. Since $X$ is compact $d_s$ is
bounded hence, by restriction, $d_g$ is bounded.

Suppose that $d_s$ is unbounded. Then for $M>0$ there exist $f\in C(X)$ and
$x,y\in X$ such that $|f(x)-f(y)|\geq M$ and $\|[D,\pi_\tau(f)]\|\leq
1$. By continuity of $f$ exist $v,w\in V$ such that
$|f(x)-f(\tau(v))|\leq 1$, $|f(y)-f(\tau(w))|\leq 1$. Hence
$|f(\tau(v))-f(\tau(w))|\geq M-2$ which implies $d_s(\tau(v),\tau(w))\geq M-2$.
It follows that $d_g$ is unbounded. So boundedness of $d_g$ implies
boundedness of $d_s$. 
\qed

\medskip
At this point, even if the graph is bounded, we neither know whether $d_s$ is continuous nor whether finite $\|f\|_{Lip,d_s}$ implies finite $\|f\|_{Lip,d} $. 
%%%%%%%%%%%
\begin{cor} $d_s$ is Lipschitz equivalent to $d$ if and only if
there exists $C>0$ such that $d_g(v,w)\leq C d(\tau(v),\tau(w))$. In this case
$\|f\|_{Lip,d}\leq C\|f\|_{Lip,d_s} $. 
\end{cor}
%%%%%%%%%%%
\noindent{\em Proof:} 
By Lemma~\ref{lem-comparison}
$d_s(x,y)  
= \sup \{|f(x)-f(y)|:  {\forall v,w\in V}\:
|f(\tau(v))-f(\tau(w))|\leq d_g(v,w)\}$
Using  $d_g(v,w)\leq C d(\tau(v),\tau(w))$ we get
%%%
\begin{multline*}
d_s(x,y)  
\leq C \sup \Big\{ |f(x)-f(y)|\ : \ \\  \forall v,w\in V , \,
|f(\tau(v))-f(\tau(w))|\leq d(\tau(v),\tau(w)) \Big\}.
\end{multline*}
%%%
The latter equals $Cd(x,y)$ by continuity of $d$ and $f$. The other statements are clear.
\qed

%%%%%%%%%%%%%%%%%%%%%%%%%%%%%%%%%%%%%%%%%%%%%%%%%%%%%%%%%%%%%%%%%%%%%%%%%%
\section{Spectral triples for compact ultra metric spaces}

We now suppose that $d$ is an ultra metric. This means that
\begin{equation}\label{ultra} 
d(x,y) \leq \max\{d(x,z),d(y,z)\} 
\end{equation}
for all $x,y,z\in X$. We denote by $B_\delta(x)=\{y\in
X:d(x,y)<\delta\}$ the open $\delta$-ball. 
%%%%%%%%
\begin{lemma} 
Let $(X,d)$ be a compact ultra metric space.
\begin{enumerate}
\item For any $\delta$, either $B_\delta(x)=B_\delta(y)$ or
 $B_\delta(x)\cap B_\delta(y)=\emptyset$. 
In particular there is a unique cover by $\delta$-balls
($\delta$-cover). Moreover this cover is a partition and hence $X$ is totally disconnected.  
\item The image of $d$ is discrete away from $0$. $0$ is an
 accumulation point, provided $X$ is infinite. 
In other words, there exists a strictly decreasing sequence $(\delta_n)_n$ converging to zero such that 
$\mbox{\rm im}\, d = \{\delta_n :n\in\N\}$. 
\end{enumerate}
\end{lemma}
%%%%%%%%
As a result of the first part of the lemma 
the number of closed $\delta$-balls needed to
cover $X$ is uniquely defined and so defines a function 
$N_X(\delta)$, somewhat like the patch-counting function in the context of discrete tiling spaces.
It is a positive, decreasing, semi-continuous, integer-valued function on $\R^{>0}$. 
\medskip

\noindent
{\em Proof:} 
Suppose $z\in B_\delta(x)\cap B_\delta(y)$.  Then (\ref{ultra})
implies $d(x,y) < \delta$. Furthermore for any $y'\in B_\delta(y)$
we have $$d(x,y')\leq\max \{d(x,y),d(y,y')\} < \delta.$$ Hence
$B_\delta(y)\subset B_\delta(x)$ and by a symmetric argument the first
statement follows.
The uniqueness of $\delta$-covers is then obvious as is the fact that any cover is a partition.
In particular, $\delta$-balls are clopen and hence, by the very definition, $X$ is totally disconnected. 

As $N_X$ is a positive, decreasing, semi-continuous, integer-valued function it has on any interval $[\delta,\infty)$, $\delta>0$, only finitely many jumps.
Furthermore, the image of $d$ coincides with its discontinuities, i.e.\ with the points at which $N_X$ jumps by at least $-1$.
Hence $[\delta,\infty)\cap \mbox{\rm im}\, d$ is finite. \qed 

%\bigskip

%%%%%%%%%%%%%%%%%%%%%%%%%%%%%%%%%%%%%%%%%%%%%%%%%%%%%
\subsection{The Michon tree and neighborhood graph} \label{sec-3.1} 

There is a tree $\Tt = (\Tt^{(0},\Tt^{(1)})$ associated with
the ultra metric compact space $(X,d)$. Its construction goes back to
Michon \cite{Mich85} and
has been already used by Pearson and Bellissard for the same purpose, 
although we use here the non-reduced version. 

Let $(\delta_n)_n$ be a strictly decreasing sequence tending to zero such that
$\im d = \{\delta_n:n\in\N\}$. We define the vertices
$\Tt^{(0)}=\bigcup_{n=0}^\infty\Tt_n^{(0)}$ of a tree using
this sequence. The set of level $n$ vertices $\Tt_n^{(0)}$ is the closed
$\delta_n$-cover, i.e.\ a vertex corresponds to a clopen
$\delta_n$-ball. We denote the clopen subset 
corresponding to $v\in \Tt_n^{(0)}$ by $[v]$. $\Tt_0^{(0)}$ contains a
single vertex, denoted $\circ$ and called the root, which stands for
all of $X$.

The edges of the tree are oriented going from vertices of level $n$ to vertices of level $n+1$: in fact there is an edge from $v \in \Tt^{(0)}_n$ to $w\in \Tt^{(0)}_{n+1}$ iff $[w]\subset[v]$.
Any vertex has thus one incoming edge. 
A branching vertex is a vertex which has at least two outgoing edges.

Consider the set of infinite rooted paths $\Pi_\infty$ on the Michon tree.
An infinite path $\xi$ is a sequence of vertices $\xi_n\in \Tt^{(0)}_n$ such that $(\xi_n,\xi_{n+1})$ forms an edge.
Equipped with the relative topology of the product topology $\Pi_\infty$ is a compact space.
The map $\xi\mapsto x(\xi)$ which assigns to $\xi$ the unique point in the infinite intersection $\bigcap_n [\xi_n]$ is a homeomorphism between $\Pi_\infty$ and $X$.

We will need what we call the neighborhood graph $H = (\Tt^{(0)},H^{(1)})$  of $(X,d)$.
It has the same vertices as the Michon tree, but its edges are ``horizontal'' if one views the tree ``vertically'', that is, between vertices of the same level.
Let $v\in \Tt^{(0)}_{n-1}$ be a branching vertex.
We denote by $\Tt^{(0)}(v)\subset \Tt^{(0)}_n$ the set of vertices which are linked to $v$ by an edge in $\Tt^{(1)}$.
We introduce an un-oriented edge between any two vertices of $\Tt^{(0)}(v)$. 
An unoriented edge counts simply as two edges oriented in opposite direction.
We denote by $H^{(1)}_n$ the collection of such edges between vertices of level $n$.

%%%%%%%%%%%%%%%%%%%%%%%%%%%%%%%%%%%%%%%%%%%%%%%%%%%%%
\subsection{Approximation graphs} 
\label{approxgraph}
Like Pearson and Bellissard, we introduce choice functions.
Our choice functions are, however, not the same as the ones used in \cite{PB09}, see Section~\ref{sec-comparison} for a detailed comparison. 
%%%%%%%%
\begin{defn}
By a collection of choices we mean a collection $(\tau_v)_{v\in \Tt^{(0)}}$ of maps $\tau_v:\{v\}\hookrightarrow \Tt^{(0)}(v)$. 
In particular $(v,\tau_v(v))$ is a vertical edge in the Michon tree. 
Note that a collection of choices defines a map $\tau:\Tt^{(0)}\to
\Pi_\infty$. $\xi=\tau(v)$ is iteratively determined by 
$\xi_{n+1} = \tau_{\xi_n}(\xi_n)$. This map satisfies
\begin{enumerate}
\item $\tau(v)$ goes through $v$. 
\item If $w\in \tau(v)$ then $\tau(w) = \tau(v)$.
\end{enumerate}
We call $\tau$ a choice function. Conversely, any choice
function gives rise to a collection of choices.
\end{defn}
%%%%%%%%
The graphs we are interested in are the graphs defined by a collection
of choices. 
Given such a collection we define a new graph $\Gamma(\tau)=(V,E)$. It
is obtained from the neighboring graph by quotienting out the  
equivalence relation generated by $v \sim \tau_v(v)$. Stated
differently, two vertices $v,w$ in $\Tt^{(0)}$ are identified if
$\tau(v)=\tau(w)$. We denote by $V$ the equivalence classes in
$\Tt^{(0)}$ and by $E$ the induced edges,  that
is, we draw an unoriented edge between two vertices of $V$ if they
have representatives in some $\Tt^{(0)}_n$ which form an edge. The edges
are thus still naturally grouped into levels, $E=\bigcup_n E_n$ where 
$E_n$ corresponds to $H^{(1)}_n$ ($E_n$ may be empty). 
Note that $\Gamma(\tau)$ is a connected graph.
We consider it as a metric graph by associating to an edge $e\in E_n$ the length $\delta_n$. 
We call $\Gamma(\tau)$ an approximation graph, because its vertices approximate the space and its edges define which points are neighbors.

If we combine a choice function with the natural homeomorphism $\xi\mapsto x(\xi)$ between $\Pi_\infty$ and $X$ we obtain a dense embedding of $\Gamma(\tau)$ in $X$. 
We obtain thus a spectral triple over $C(X)$ which we call the spectral triple for $\Gamma(\tau)$ as described in Section~\ref{embdgraph}.

%%%%%%%%%%%%%%%%%%%%%%%%%%%
\subsubsection{Example: aperiodic one-sided subshifts}
Recall the definition of a one-sided subshift from Section~\ref{tilings}.
We may associate to it the tree of words defined as follows:  the vertices of level $n$ are the words of length $n$.
There is a vertical edge between a word $w$ (possibly the empty one) and any of its one-letter extensions $wa$, $a\in\A$.
In particular, a branching vertex is a word which can be extended by one letter in more than one way. Such a word is called right special. 
Note that aperiodicity implies that there is at least one right special word per length (see the proof of Lemma~\ref{lem-repulsive}).

Recall that the subshift space $\Xi_\L$ has a topology 
which can be generated by an ultra metric of the form
%%%%%%%%
\begin{equation}
\label{ssmetric}
d(\xi,\eta) = \inf\{\delta_n \, : \, \forall m\leq n : \xi_m=\eta_m\}\,,
\end{equation}
%%%%%%%%
where $(\delta_n)_n$ is a strictly decreasing sequence tending to $0$. 
The Michon tree associated to $(X,d)$ is simply the tree of words. 
It does actually not depend on the choice of the sequence $(\delta_n)_n$ as long as the latter is strictly decreasing and tending to $0$. 
Clearly the Michon tree has uniformly bounded branching.

We describe the neighboring graph $H$.
There is a horizontal edge between two distinct words of equal length $n$ whenever their prefices of length $n-1$ coincides.
In particular, horizontal edges appear precisely between different extensions of right-special words. 

To extend a right special word one needs to choose. This is exactly what a collection of choices
$(\tau_v)_v$ is good for, it choses for each word a one-letter extension.  A
choice function associates thus to any word an infinite extension.   
$\Gamma(\tau)$ is therefore a graph whose vertices yield a dense set
of points in $\Xi$ corresponding to the chosen infinite extensions of
words (many of which have the same extension).
An edge of   $\Gamma(\tau)$ of level $n$ appears between two such points which are extensions of a right special word of length $n-1$ whose $n$-th letter differs.
The larger the $n$ the smaller their distance. 

%%%%%%%%%%%%%%%%%%%%%%%%%%%%%%%
\subsubsection{Example: tilings of finite local complexity}
Consider a tiling of finite local complexity.
It defines a sequence of sizes $(r_n)_n$.
We consider the tree of patches: its level $n$ vertices are the $r_n$-patches and its root represents the empty patch ($r_0=0$). 
We introduce a vertical edge between an $r_n$-patch and any
of its extensions to an $r_{n+1}$-patch. Given a metric as in
Proposition~\ref{prop-metric} defined by a function $\delta(r)$
the Michon tree associated to $(X,d)$ coincides with the above tree,
except in the case that, for some $n$, 
all $r_n$ patches force their extensions to
$r_{n+1}$-patches. In the latter case Michon tree is obtained from the tree of patches by contracting the vertices of levels which do not have branching vertices. 
Note that if the tiling has a uniform bound on the number of possible patch extensions, then the Michon tree has uniformly bounded branching.

As for the neighboring graph, the horizontal edges
correspond exactly to the pairs of distinct $r_{n}$-patches which
are extensions of the same   $r_{n-1}$-patch. Finally, the choice function
associates to each $r_n$-patch a tiling of which it is the center.

%%%%%%%%%%%%%%%%%%%%%%%%%%%%%%%%%%%%%%%%%%%%%%%%%%%%%%%%%%%%%%%%%%%%%%%%%%%%
\section{The spectral distance and aperiodic order}
From now on, we consider only ultra metric spaces $(X,d)$.

%%%%%%%%%%%%%%%%%%%%%%%%%%%%%%%%%%%%%%%%%%%%%%%%%%%%%
\subsection{Spectral distance}
We now aim at a formula for $d_s$.
For that we use the description of $X$ as the set $\Pi_\infty$ on infinite rooted paths in the Michon tree.
Recall that a path $\xi\in\Pi_\infty$ is a sequence $(\xi_n)_n$ of vertices where $\xi_n\in \Tt^{(0)}_n$.
Given a path $\xi \in\Pi_\infty$ and integers $n<m$, we write $\xi_{[n,m]}$ for its restriction from level $n$ up to $m$.
Note that $\xi = \tau(v)$, $v\in \Tt^{(0)}_{n_0}$ iff  
$\xi = \tau(\xi_n)$
for all $n\geq n_0$. 
Now define $\b^\tau_n:\Pi_\infty\to \{0,1\}$ by 
$$\b^\tau_n(\xi) = \left\{\begin{array}{ll}
0 & \mbox{if }\xi_{n+1} = \tau(\xi_n)_{n+1}\\
1 & \mbox{otherwise}\end{array}\right. .$$
Two distinct paths $\xi$ and $\eta$ will become distinct at a certain vertex. We denote this vertex by
$\xi\vee\eta$, it is thus the vertex of highest level which is common to $\xi$ and $\eta$.
We furthermore denote by $|\xi\vee\eta|$ the level of this vertex.
%%%%%%%%%%
\begin{lemma}
Consider the Michon tree for an ultra metric compact space. 
We denote by $(\delta_n^\tau)_n$ the decreasing sequence of the image of the metric.
Let $d_s^\tau$ be the spectral distance of the triple defined by that metric and a chosen
choice function $\tau$. If $d_s^\tau$ is continuous then
$$ d_s^\tau(\xi,\eta) = 
%d(\xi,\eta) 
\delta_{|\xi\vee\eta|} 
+ \sum_{n>|\xi\vee\eta|} \b^\tau_n(\xi) \delta_n 
+ \sum_{n>|\xi\vee\eta|}\b^\tau_n(\eta) \delta_n.$$ 
\end{lemma}
%%%%%%%%%%
{\em Proof:} Under the hypothesis of continuity of $d_s$ it suffices to
show that the formula is correct provided $\xi,\eta\in\im \tau$. 
%In that case $d_s(\xi,\eta) = d_g(\xi,\eta)$.
Suppose that $\xi\vee\eta$ has level $n-1$ (so $(\xi_n,\eta_n)$ is a horizontal edge).
Let $h_1\cdots h_N$ be a shortest path on $\Gamma(\tau)$ linking $\xi$ to $\eta$. 

We claim that $h=h_1\cdots h_N$ passes through the edge $(\xi_n,\eta_n)$. To see this first note that 
if $h$ contains an edge from $E_m$, $m<n$, then it contains a cycle and hence would not be of shortest length. Second, if $h$ does not pass through $(\xi_n,\eta_n)$ then, in order to get from $\xi_n$ to 
$\eta_n$, which it certainly must, it must contain a sub-path from $\xi_n$ to 
$\eta_n$ which passes through at least two other edges of $E_n$. By replacing
that sub-path by the edge $(\xi_n,\eta_n)$ $h$ could thus be made shorter. 

The above shows that 
%%%
\[
d_g(\xi,\eta) = d_g(\xi,\tau(\xi_n)) + d_g (\tau(\xi_n),\tau(\eta_n))+ d_g (\tau(\eta_n),\eta_n) \,.
\]
%%%
Clearly $d_g (\tau(\xi_n),\tau(\eta_n)) = l((\xi_n,\eta_n))= \delta_{\xi\vee\eta}$.
Furthermore, $|\xi\vee\tau(\xi_n)|\geq n$ and so we may recursively apply the above argument to obtain
$d_g(\xi,\tau(\xi_n)) =  \sum_{n>|\xi\vee\eta|} \b^\tau_n(\xi) \delta_n$.\qed
\bigskip

%%%%%%%%%%
\begin{cor}
\label{cor-cont}
If $\sum_n\delta_n<\infty$ then $d_s^\tau$ is continuous.
\end{cor}
%%%%%%%%%%
{\em Proof:} For fixed $n$, $\xi\mapsto \b^\tau_n(\xi)$ is continuous since it depends only on 
$\xi_{[0,n+1]}$. Under the above condition \(\sum_n \sup_\tau(\b^\tau_n) \delta_n \leq \sum_n\delta_n<\infty\) showing that the series is normally convergent. 
Hence $\xi \mapsto \sum_n\b^\tau_n(\xi)\delta_n$ is continuous.\qed
\bigskip

It seems to us more interesting, not to suppose that the series $\sum_n\delta_n$ is summable.
We will see below that the question of summability of that series
gives rise to an interesting characterization of Sturmian sequences,
see Theorem~\ref{thm-cexstur}. 
\begin{cor}
A necessary condition for $d_s^\tau$ to be continuous is that there exists a $C>0$ such that for all
$\xi\in\Pi_\infty$ 
$$\sum_n\b^\tau_n(\xi)\delta_n\leq C.$$ 
\end{cor}
{\em Proof:} If $d_s$ is continuous and $\eta\neq \xi$ then the above
formula tells us that 
$\sum_{n>|\xi\vee\eta|}\b^\tau_n(\xi)\delta_n\leq d_s(\xi,\eta)$. 
Furthermore, by compactness of $X$ the latter is bounded.
\qed

%Note that $\sum_n\b^\tau_n(\xi)\delta_n\leq\zeta^{low}(1)$ but the r.h.s.\ may be infinite.
\begin{cor} \label{cor-equiv}
$d_s^\tau$ is Lipschitz equivalent to $d$ if and only if
there exists a constant $C>0$ such that, for all $\xi\in\Pi_\infty$ and all $n_0$ for which $\xi_{n_0}$ is branching we have
$$\sum_{n> n_0}
\b^\tau_n(\xi) \delta_n \leq C d(\xi,\tau(\xi_{n_0})).$$ 
\end{cor}
{\em Proof:} We have
\begin{eqnarray*}
d_s(\xi,\eta) & = & \delta_{\xi\vee\eta} +   
\sum_{n>|\xi\vee\eta|}
\b^\tau_n(\xi) \delta_n + \sum_{n>|\xi\vee\eta|}
\b^\tau_n(\eta) \delta_n
%\leq (1+2C) \delta_{\xi\vee\eta} 
\\ &\leq & d(\xi,\eta) + C\left(d(\xi,\tau(|\xi\vee\eta|)) + d(\eta,\tau(|\xi\vee\eta|))\right)
\end{eqnarray*} 
and the corollary follows from $ d(\xi,\tau(|\xi\vee\eta|)) \leq d(\xi,\eta)$.
\qed
\medskip

It is clear the the spectral metric depends on $\tau$. Let us consider therefore the extremes
%%%
\[
 \underline{d}_s := \inf_{\tau\in Y} d_s^\tau , \qquad  \overline{d}_s := \sup_{\tau\in Y} d_s^\tau \,.
\]
%%%
which depend canonically on $(X,d)$.
The question as to whether, given a compact ultra metric space $(X,d)$, $\overline{d}_s$ is continuous or even Lipschitz-equivalent to $\underline{d}_s$ can be seen as a characterization of compact ultra metric spaces.
%%%%%%%
\begin{cor}\label{cor-inf}
$$\inf_{\tau\in Y} d_s^\tau(\xi,\eta) = d(\xi,\eta).$$
\end{cor}
%%%%%%%
{\em Proof:} 
If $\tau$ is a choice function such that $\tau(\xi_{|\xi\vee\eta|+1})=\xi$ and
$\tau(\eta_{|\xi\vee\eta|+1})=\eta$ then the sums in the r.h.s.\ of the formula for $d_s(\xi,\eta)$ are nil.
\qed
\medskip 

Define $a(\xi_n) = \# \Tt^{(0)}(\xi_n) - 1$, i.e.\ the branching number of the vertex $\xi_n$ {\bf minus} $1$, and set
\(\overline{\b}_n(\xi) = \sup_\tau \b_n^\tau(\xi)\)
so that \(\overline{\b}_n(\xi) =1\) if $a(\xi_n)>0$ and $0$ else.
In particular, $\overline{\b}_n(\xi)$ is $1$ whenever the vertex $\xi_n$ branches.
%%%%%%%%
\begin{lemma}
\label{lem-choice}
Given any sequence $(c_n)_n\subset \{0,1\}$ such that $c_n \leq \overline{\b}_n(\xi)$, there exists
a choice function $\tau$ such that $c_n = \b^\tau_n(\xi)$.
\end{lemma}
%%%%%%%%
{\em Proof:} 
Any choice function satisfying $\tau(\xi_n)_{n+1}\neq \xi_{n+1}$ provided $c_n=1$ does the job.
\qed
\medskip

This allows us to rewrite the above criteria as follows:
\begin{cor}\label{cor-Lip}${}$
\begin{enumerate}
\item
If $\sup_\xi\sum_n \overline{\b}_n(\xi)\delta_n=\infty$ then $\overline{d}_s$ is not continuous. 
\item $\overline{d}_s$ is Lipschitz equivalent to $d$ if and only if there exists a constant $C>0$ such that, for all $\xi\in\Pi_\infty$ and $m$ 
$$\sum_{n> m}
\overline{\b}_n(\xi) \delta_n \leq C \delta_{m}.$$ 
\end{enumerate}
\end{cor}

%%%%%%%%%%%%%%%%%%%%%%%%%%%%%%%%%%%%%%%%%%%%%%%%%%%%%%%%%%%%%
\subsection{Episturmian subshifts}
We start with a simple result which shows that the exponential choice for the metric is not useful for our purpose.
%%%%%%%%%
\begin{thm}
\label{thm-exp}
Given any one-sided subshift. If $\delta_n=e^{-n}$ then $\bar{d}_s$ is Lipschitz equivalent to $d$.
\end{thm}
%%%%%%%%%
{\em Proof:}
$$\sum_{n > m} \overline{\b}_n(\xi) \delta_n \leq \sum_{n > m} e^{-n} = \frac{e^{-m}}{e-1} = C\delta_m$$
where $C=\frac{1}{e-1}$.
\qed
\medskip

A simple class of subshifts consists of  repetitive subshifts which have the particular property that, for any $n\geq 1$, there is
exactly one right special word. We look here at the one-sided version
only. The Michon tree has then the property that there is exactly one
branching vertex per level.
Repetivity implies that the set of right special words determine the
subshift completely as the language $\L$ must be the set of sub-words
of right special words. Given that the suffix of a right special word
must be right special we can define a {\em left} infinite word $\Rr$
as the sequence whose $n$-suffix is precisely the right special word
of length $n$. $\Rr$ is called the characteristic word of the
subshift. 

Examples of subshifts of this type are episturmian sequences (for these one requires also to have precisely one left-special word per length) \cite{GJ09}, among which one finds the Arnoux-Rauzy subshifts \cite{AR91}.
Again a sub-class of the latter are 
canonical projection method tilings of the following type:
%of dimension $1$ and codimension $n$: 
A line $E$ is placed in $\R^{N}$ in such a way that it does not
intersect the integral lattice
$\Z^{N}\subset\R^{N}$. In particular, w.r.t.\ a lattice base
$E$ is parallel to a vector $\vec\nu$ whose components are
rationally independent and positive.
The subshift is characterized by $\vec\nu$.
Let $C=[0,1]^{N}$ be the unit cube in $\R^{N}$.
The set $(E+C)\cap \Z^{N}$ is then orthogonally projected onto
$\R\vec\nu$. This yields a sequence of intervals of different lengths which we
encode by the $N$ letters of an alphabet $\A$. By declaring the 
interval whose interior or left boundary contains $0\in  \R\vec\nu$ to
be $a_0$ we obtain  
a sequence $(a_n)_n\in\A^\Z$. 
The right half (in the direction of $\vec\nu$) 
of the sequence so obtained corresponds
to a one-sided sequence of the subshift. 
The subshift space characterized by $\vec\nu$
is the set of all sequences obtained in this way.
The characteristic word $\Rr$
is the left half of the sequence obtained when one uses $\R\vec\nu$ in
place of $E$. 

Again a particular sub-class of the above are the Sturmian subshifts which correspond to the case $N=2$. Thus they are determined by $\vec{\nu}=(1,\theta)$ for some irrational $\theta$ which may be taken in the interval $(0,1)$.

For Sturmian subshifts the  notions of repulsiveness are equivalent to linear recurrence, or linear repetivity \cite{DHS99,Du00}.
%%%%%%%%%%%%%%%%
\begin{lemma}
\label{lem-sturLR}
For a Sturmian subshift  associated with the irrational \(\theta \in (0,1) \), the following are equivalent.
%%%%%%%%%%%%%%%%
\begin{enumerate}[(i)]

\item Linear recurrence;

\item Repulsiveness;

%\item Right-repulsiveness;

\item $\theta$ has bounded continuous fraction expansion.

\end{enumerate}
%%%%%%%%%%%%%%%%
\end{lemma}
%%%%%%%%%%%%%%%%
{\em Proof:}
For a Sturmian subshift the characteristic word $\Rr$ can be obtained as follows:
Let $[1+\mu_0,\mu_1,\cdots]$ be the continued fraction expansion of
$\theta$. 
Let $\Aa = \{ a, b\}$, and define $\sigma_0$ and $\sigma_1$ as the substitutions 
\begin{eqnarray}
\sigma_0(a) = a\:\:&\quad& \sigma_0(b) = ba \\
\sigma_1(a) = ab&\quad& \sigma_1(b) = b
\end{eqnarray}
Then $\Rr$ is the left infinite sequence whose suffices are given by
%%%
\[
\Rr_k:= \sigma_0^{\mu_0}\sigma_1^{\mu_1} \cdots
\sigma_0^{\mu_{2k}}(b)\,, \quad k \ge 0\,.
\]
%%% 
In fact, since $\sigma_1^l\sigma_0^m(b)= b(ab^l)^m$, $\Rr_k$ is a suffix of
$\Rr_{k+1}$.
% and so we may define $\Rr=\lim_k\Rr_k$.  
%The right special words of $\L_\theta$ are precisely the suffices of
%the $\Rr$. 
Iteration shows that $\Rr$ has a suffix of the form $u_k^{\mu_k}$ where $u_k$ depends only on $\mu_0,\cdots,\mu_{k-1}$. In
particular, $\Rr_{(\mu_k-1)|u_k|} = u_k^{\mu_k-1}$ and $\ell_{(\mu_k-1)|u_k|}=
\frac{|u_k|}{ (\mu_k-1)|u_k|} = \frac{1}{\mu_k-1}$.
In particular if $\inf_n\ell_n>0$ than the continued fraction expansion is bounded, and this is a well-know characterization of linear recurrence \cite{DHS99,Du00}.
Thus we can conclude as follows:
Linear repetivity implies repulsiveness (Theorem~\ref{thm-Lenz}) implies right special repulsiveness
(Lemma~\ref{lem-repulsive})  implies bounded continued fraction expansion implies linear repetivity. Thus all four notions are equivalent. 
\qed

%%%%%%%%%
\begin{prop}
\label{prop-epistur}
Consider %the subshift space $\Xi$ of 
a repetitive one-sided subshift
which has exactly one right special word per length.
Let  $(\delta_n)_n$ be a (strictly decreasing and converging to zero) sub-multiplicative sequence, 
i.e.\ $\delta_{nm}\leq \overline{c}\delta_n\delta_m$ for some $\overline{c}>0$. 
We provide the subshift space with the metric $d$ defined by $(\delta_n)_n$ as in equation \eqref{ssmetric}.
If the subshift is repulsive then $\bar{d}_s$ is Lipschitz equivalent to $d$.
\end{prop}
%%%%%%%%%
{\em Proof:}
Consider first the case $\overline{c}=1$.
Let $\xi$ be an infinite path on the Michon tree and $n_k$ be the level of its $k$-th branching vertex. 
In other words the $n_k$ are precisely the values for which $\overline{\b}_{n_k}(\xi)=1$.
Again equivalently, the $n_k$ are precisely the values for which  $\xi_{[0,n_k]}$ is a right special word. 
It follows that $\xi_{[0,n_k]}$ is not only a prefix but also a suffix of $\xi_{[0,n_{k+1}]}$, by uniqueness of the right special words of length $n_{k+1}$.
Hence $n_{k+1}-n_k =\ell_{n_k}n_k \geq \ell_{\Rr} n_k$ where $\ell_{\Rr}=
\inf_n\ell_n$.
%\frac{rt(\omega_n)}{n}$. 
It follows that $n_k\geq (\ell_{\Rr}+1)^{k-l}n_l$. 
Thus 
%%%
\[
\sum_{n > n_l} \overline{\b}_n(\xi) \delta_n = \sum_{k>l} \delta_{n_k}
\leq \delta_{n_l} \sum_{k\geq 1} \delta_{[(\ell_\Rr+1)^k]}
\]
%%%
where $[(\ell_\Rr+1)^k]$ is the integer part of $(\ell_\Rr+1)^k$. 
Now if the subshift is right-special repulsive there exist a $k_0$
such that $(\ell_\Rr+1)^{k_0}\geq 2$. 
It follows that $ \delta_{[(\ell_\Rr+1)^{nk_0+m}]}\leq  {\delta_{2}}^n
\delta_{[(\ell_\Rr+1)^m]}$. 
Since $(\delta_n)_n$ is strictly decreasing we have $\delta_2<1$  so that    
the series is convergent and yields
the finite constant for the bound in Corollary~\ref{cor-Lip}(2)
(it suffices to test the criteria for $m = n_l$ as the contribution
for non branching vertices to the left hand side is nil).
Now if $\overline{c}>1$ we can rescale the $\delta_n$ so as to obtain
$\overline{c}=1$. This will 
change the metric $d$ by a constant factor $\overline{c}^{-\frac12}$ but not its
Lipschitz-class.\qed 

%%%%%%%%%%%
\begin{prop}
\label{prop-epistur-reverse}
Consider %the subshift space $\Xi$ of 
a repetitive one-sided subshift
which has exactly one right special word per length.
Let  $(\delta_n)_n$ be a (strictly decreasing and converging to zero) sequence satisfying 
$\delta_{2n}\geq \underline{c}\delta_n$ for some $\underline{c}>0$. 
We provide the subshift space with the metric $d$ defined by $(\delta_n)_n$ as in equation \eqref{ssmetric}.
If the subshift is {\em not} repulsive then $\bar{d}_s$ is {\em not}
Lipschitz equivalent to $d$.
\end{prop}
%%%%%%%%%%%
{\em Proof:}
We use the terminology of the proof of Prop.~\ref{prop-epistur}. 
Suppose that $\ell_{n_k}<1$. Then $n_{k+1}<2n_k$ and so the prefix and the suffix of
length $n_k$ of $\xi_{[0,n_{k+1}]}$ overlap. Let $u$ be the prefix of
$\xi_{[0,n_{k}]}$ of length $n_{k+1}-n_k$. The
overlap forces   $\xi_{[0,n_{k}]}$ to be of the form
$\xi_{[0,n_{k}]}=u^{p_k} v$ where $v$ is a possibly empty word of length
strictly smaller then $n_{k+1}-n_k$. Moreover, $v$ must be a prefix of
$u$. Hence $\xi_{[0,n_{k}]}=v(wv)^{p_k}$ for some non-empty word $w$.
It follows that $p_k$ is the integer part of $\ell^{-1}_{n_k}$ and that
$\xi_{[0,n_{k+1}]}=v(wv)^{p_k+1}$. But this implies that
$\xi_{[0,|v|+j|u|]}$ are right special words for all $j=0,\cdots
p_k+1$. Let $q_k$ be the integer part of $\frac{p_k}2$ and
$m_k=q_k(n_{k+1}-n_k)$. Then 
$$
\sum_{n\geq m_k} \overline{\b}_n(\xi) \delta_n \geq  
\sum_{j=q_k}^{p_k}\delta_{|v|+j|u|}
\geq  \sum_{j=q_k+1}^{2q_k}\delta_{j|u|}.
$$
Hence
\[
\delta^{-1}_{m_k} \sum_{n\geq m_k} \overline{\b}_n(\xi) \delta_n \geq
\sum_{j=q_k+1}^{2q_k}\frac{\delta_{j|u|}}{\delta_{q_k|u|}}
\geq
\sum_{j=q_k+1}^{2q_k}\frac{\delta_{2q_k|u|}}{\delta_{q_k|u|}}\geq
\underline{c} q_k.
\]
Since $\inf_n \ell_n = 0$ if and only if $\sup_k q_k = +\infty$
the statement follows from Corollary~\ref{cor-Lip}(2).
\qed

\medskip
The previous two propositions combined yield our main result:
%%%%%%%%%%%
\begin{thm}
\label{thm-epistur}
Consider a repetitive one-sided subshift which has exactly one right special word per length. 
Let  $(\delta_n)_n$ be a (strictly decreasing and converging to zero) sequence such that 
$\underline{c}\delta_n\leq\delta_{2n}$ and  $\delta_{nm}\leq \overline{c}\delta_n\delta_m$ 
for some $\overline{c},\underline{c}>0$. 
We provide the subshift space with the metric $d$ defined by $(\delta_n)_n$ as in equation \eqref{ssmetric}. 
The following are equivalent:
\begin{enumerate}
\item\label{item-rL}   The subshift is repulsive.
\item\label{item-Lr}   $\overline{d}_s$ is Lipschitz equivalent to $d$.
\end{enumerate}
\end{thm}
%%%%%%%%%%%%%%%

With Lemma~\ref{lem-repulsive}, in the particular case of Sturmian subshifts, we get the following.
%%%%%%%%%%%%%%%
\begin{cor} 
\label{cor-stur}
Consider a Sturmian subshift associated with the irrational \(\theta
\in (0,1)\) and a metric on its subshift space as in Theorem~\ref{thm-epistur}. 
The following are equivalent:
%%%%%%%%%%%%%%
\begin{enumerate}[(i)]
\item The Sturmian sequence is linearly recurrent;
\item The continued fraction expansion of $\theta$ is bounded;
\item $\bar{d}_s$ is Lipschitz equivalent to $d$.
\end{enumerate}
%%%%%%%%%%%%%%
\end{cor}
%%%%%%%%%%%%%%%
Note that $\delta_n=e^{-n}$ does not satisfy the hypothesis of the
theorem (as it should be). Possible choices for the $\delta_n$ satisfying
the hypothesis are
$$\delta_n = \frac{\ln^bn}{n^a}$$ for $a>0$, $b\geq 0$.
\bigskip

Recall from Lemma~\ref{cor-cont} that if $\sum_n\delta_n$ is summable 
then $\overline{d}_s$ is always continuous.
%%%%%%%%%%%
\begin{thm}
\label{thm-cexstur}
Consider a strictly decreasing sequence $(\delta_n)_n$ such that
$\sum_n\delta_n$ is not summable. 
There exist Sturmian subshifts for which
the metric  $\overline{d}_s$ is not continuous.
For these subshifts there exists even a choice function $\tau$
for which $d^\tau_s$ is not continuous. 
\end{thm}
%%%%%%%%%%%
{\em Proof:}
The same calculation as above shows that
%%%
\[
\sum_{n> 0} \overline{\b}_n(\xi) \delta_n \geq  
\sum_{j=1}^{p_k+1}\delta_{j|u|}. 
\]
%%%
Note that, since $\delta_n$ is a strictly decreasing sequence,
all series  $\sum_n\delta_{ln}$ for $l\in\N$ diverge to infinity as well.
Consider a Sturmian subshift.
Given $i$ let $\xi^{(i)}$ be an infinite path with prefix $u_i^{\mu_i}$ (notation as in the proof of Lemma~\ref{lem-sturLR}). 
There is a $n_k$ which is of the form $n_k = (\mu_i-1)|u_i|$ and
$n_{k+1} = \mu_i|u_i|$. Thus 
\begin{equation}
\sum_{n> 0} \overline{\b}_n(\xi^{(i)}) \delta_n \geq  
\sum_{j=1}^{\mu_i}\delta_{j|u_i|}. 
\end{equation}
Since $u_i$ depends only on $\mu_0,\cdots,\mu_{i-1}$ we may construct
a sequence $(\mu_i)_i$ such that the r.h.s.\ above
diverges for $i\to\infty$. 
It follows that, for a subshift defined by such a sequence, 
$\sup_i\sum_{n> 0} \overline{\b}_n(\xi^{(i)}) \delta_n=\infty$
implying by  Corollary~\ref{cor-Lip}(1) that  
$\overline{d}_s$ is not continuous.

By compactness of $\Xi$ the sequence $(\xi^{(i)})_i$ admits a limit
point $\xi$.  
By definition of the topology this point $\xi$ has arbitrarily long suffices of the type $|u_i|^{\mu_i}$.
Hence  $\sum_{n> 0} \overline{\b}_n(\xi) \delta_n$ diverges. By Lemma~\ref{lem-choice} there is a choice function $\tau$ such that $ \overline{\b}_n(\xi)={\b}^\tau_n(\xi) $. \qed

%%%%%%%%%%%%%%%%%%%%%%%%%%%%%%%%%%%%%%%%%%%%%%%%%%%%%%%%%
\subsection{FLC-tilings with equidistributed patch frequencies} 

Consider a tiling of finite local complexity which we view as an infinite path $\xi$ in the Michon tree (or the tree of patches).
We denote $\mbox{freq}_n^{min}$, $\mbox{freq}_n^{max}$ the minimal and maximal frequencies of
$r_n$-patches.
If the tiling has equidistributed frequencies then the frequencies of $r$-patches are of the same asymptotic behavior and one is given by the inverse of the patch counting function.
We denote by $p_n(\xi)$ the $r_n$-patch associated with the vertex $\xi_n$, so that $a(\xi_n)+1$ is the  number of $r_{n+1}$-patches extending $p_n(\xi)$.
To obtain an estimate for the number of non-zero $a(p_n(\xi))$ we consider the equation
$$\mbox{freq}(p_n(\xi)) = \sum_{\tilde p} \mbox{freq}(\tilde p) =
\mbox{freq}(p_{n+1}(\xi)) + \sum_{\tilde p\neq p_{n+1}(\xi)}
\mbox{freq}(\tilde p) $$
the sum running over all $r_{n+1}$-patches $\tilde p$ extending $p_n(\xi)$. Using 
$$\sum_{\tilde p\neq p_{n+1}(\xi)}
\mbox{freq}(\tilde p)\geq  a(\xi_n) \mbox{freq}_{n+1}^{min}$$
we get
$$\mbox{freq}_{n}^{max}
\geq 
\mbox{freq}(p_n(\xi)) \geq \sum_{k\geq n} a(\xi_k)\mbox{freq}_{k+1}^{min}.$$
Let $(N_j)_j$ be a strictly increasing sequence of natural numbers. Then the above yields
$$\mbox{freq}_{N_j}^{max}
\geq \sum_{k = N_j}^{N_{j+1}-1} a(\xi_k)\mbox{freq}_{N_{j+1}}^{min}.$$
The following results exploit this formula. We suppose that the patch counting function 
satisfies
%sub-multiplicative if, for each $a$ there is a constant $\tilde{c}(a)$ such that
%%%%%%%
\begin{equation}\label{eq-sub}
P(4r)\leq \tilde{c}P(r)
\end{equation} 
%%%%%%%
for some constant $\tilde c$ and sufficiently large $r$.
For example polynomial patch counting functions satisfy this bound asymptotically.
In particular, linearly repetitive tilings have equidistributed frequencies and have patch counting function which satisfies (\ref{eq-sub}). 
%%%%%%%%%%
\begin{lemma}
\label{lem-equidist-UBPE}
Consider a tiling which has equidistributed frequencies and whose patch counting function satisfies (\ref{eq-sub}).
Then the tiling has a uniform bound of its number of possible patch extensions, and thus the tree of patches as uniformly bounded branching.
\end{lemma}
%%%%%%%%%%
{\em Proof:} 
%Asymtotically, $\mbox{freq}_{n}^{max}$ and $\mbox{freq}_{n}^{min}$ behave like $R^{-d}$.
If we take $N_j$ such that $2^j\leq r_{N_j}  = 2^{j+1}$. This is possible for large enough $j$ as $r_{n+1}-r_n$ is bounded by finite local complexity. Then the
ratio $\mbox{freq}_{N_j}^{max} / \mbox{freq}_{N_{j+1}}^{min} $ is
bounded by $\frac{C}{c} \frac{p(r_{N_{j+1}})}{p(r_{N_j})}\leq
\frac{C}{c} \tilde c$.  
Hence
$$ \sum_{k = N_j}^{N_{j+1}-1} a(\xi_k)  \leq  \frac{C}{c} \tilde c.$$
Since $\xi$ was arbitrary this implies $a(\xi_n) \leq  \frac{C}{c} \tilde c$ for all $\xi$ and $n$.
\qed
\begin{thm}
\label{thm-Lip}
Consider a tiling which has equidistributed frequencies and whose patch counting function satisfies (\ref{eq-sub}). 
Suppose that the function $\delta$ which defines the metric as in Prop.~\ref{prop-metric} lies in
$L^1([1,\infty),\frac{dx}{x})$ and is submultiplicative: $\delta(ar)\leq \delta(a)\delta(r)$ for sufficiently large $a$ and 
$r$. Then $\bar d_s$ and  $\underbar d_s$ are Lipschitz-equivalent.
%distance defined by $\delta$. 
\end{thm}
{\em Proof:} We verify the criterium of Corollary~\ref{cor-equiv}.
Let $\xi\in\Pi_\infty$ and $n_0$ such that $\xi_{n_0}$ is a branching vertex.  Choose a sequence $N_j$ such that $2^j\leq \frac{r_{N_j}}{r_{n_0}}\leq 2^{j+1}$.
As in the last proof we obtain
$\mbox{freq}_{N_j}^{max} / \mbox{freq}_{N_{j+1}}^{min} \leq \frac{C}{c}\tilde c$. Since $\b_n^\tau(\xi)\leq a(\xi_n)$  we obtain  
$$ \sum_{n> n_0}\b^\tau_n(\xi) \delta_n \leq \sum_{j\geq 0} \sum_{k = N_j}^{N_{j+1}-1} a(\eta_k)\delta(r_k) \leq  \sum_{j} \delta(2^j r_{k_0}) \frac{C}{c} \tilde c \leq \tilde C \delta_{k_0}$$
where
$$ \tilde C = \frac{C}{c}\tilde c \sum_{j\geq 0} \delta(2^j).$$
Hence the Lipschitz equivalence is guaranteed
if the series is summable. This is equivalent to the integrability of
$\delta$ at infinity for the measure $\frac{dx}{x}$.\qed
\bigskip

Note that the above theorem covers any metric defined by a function $\delta$ which behaves asymptotically as $r^{-\alpha}$, $\alpha>0$.

%%%%%%%%%%%%%%%%%%%%%%%%%%%%%%%%%%%%%%%%%%%%%%%%%%%%%%%%%%%%%%%%%%%%%%%%%
\section{Zeta-functions and Laplacians for compact ultra metric spaces}

%%%%%%%%%%%%%%%%%%%%%%%%%%%%%%%%%%%%%%%%%%%%%%%%%%%%%
\subsection{Zeta-functions}
We consider now the zeta-function $\zeta(s) = \sum_{e \in E} \ell(e)^s$ and in particular its abscissa of convergence $s_0$. 
We clearly have
%%%
\[
\zeta(s) = \sum_n \sum_{e \in E_n} \ell(e)^s = 
\sum_n \sum_{v\in \Tt^{(0)}_n} 
%\left({a(v)+1} \atop 2 \right) 
a(v)(a(v)+1)\delta_n^{s} \,,
\]
%%%
where \( a(v) = \# \Tt^{(0)}(v) -1\) is the branching number of $v$ minus 1.
Using $2 a(v)\leq a(v)(a(v)+1)$, we obtain as lower bound the function 
%\left({a(v)+1} \atop 2 \right)$ 
%and consider the function
\[
\zeta^{low}(s) := 2\sum_n \s(n) \delta_n^{s}
\]
where 
\[
\s(n) := \sum_{v\in \Tt^{(0)}_n} a(v) = P(r_{n+1})- P(r_n).
\]
In particular $\zeta^{low}(s)$ is finite iff the sequence
$(\delta_n)_n$ lies in the weighted $\ell^s$ space $\ell^s(\s)$. Hence
the abscissa of convergene $s_0^{low}$ of $\zeta^{low}(s)$ is
given by
$$s_0 = \inf\{s : (\delta_n)_n\in \ell^s(\s)\}.$$
%%%%%%%%%%%%%%
We say that the Michon tree has uniformly bounded branching if there exists a constant $B$ such that $a(v)\leq B$ for all vertices $v$.
\begin{lemma} If the Michon tree has uniformly bounded branching then 
the abscissa of convergence for $\zeta$ coincides with that of $\zeta^{low}$.
\end{lemma}
%%%%%%%%%%%%%%
{\em Proof:}
Since $a(v)\leq B$ we have %$\left({a(v)+1} \atop 2 \right) \leq \frac{b+1}{2} a(v)$ and thus
$\zeta^{low}\leq \zeta \leq ({B+1} )\zeta^{low}$.
\qed

\medskip
For one-dimensional subshifts $a(v)+1$ is always bounded by the number of elements in the alphabet.  

We wish to say more about $s_0^{low}$ and $s_0$.
We start with the simple observation that if $f:\R^+\to\R^+$ is an integrable piecewise $C^1$-function such that $f'(r)\leq 0$ then there exists  $c>0$ such that $f(x)\leq\frac{c}{x}$. 
Indeed, 
$$xf(x) - f(1) = \int_1^x (f(t) + tf'(t))dt \leq \|f\|_1$$
as $ tf'(t)\leq 0$. Likewise, if $(\delta_n)_n$ is a decreasing
sequence which is summable then 
there exists a $c>0$ such that $\delta_n\leq\frac{c}{n}$.

We let
$$ \underline{\eta} =
\sup\{\gamma:\left(g(n)^{\frac{-1}{\gamma-1}}\right)_n  \mbox{
  is summable}\}$$ 
and
$$\overline{\eta} = \inf\{\gamma :
\left(\frac{g(n)}{n^{\gamma-1}}\right)_n  \mbox{ is bounded}\}$$ 
%%%%%%%%%%%%%%
\begin{thm}
\label{thm-zeta-1} Consider %the function $\zeta^{low}$ of 
a $d$-dimensional tiling of finite local complexity. 
Suppose that  $(\delta_n)_n\in \ell^{\alpha+\epsilon}\backslash
\ell^{\alpha-\epsilon}$ for some $\alpha>0$, for all $\epsilon>0$ (sufficiently small).  Then
$$ \alpha\underline{\eta} \leq s_0^{low}\leq  \alpha\overline{\eta}.$$
\end{thm} 
%%%%%%%%%%%%%%
{\em Proof:} It suffices to consider the case $\alpha=1$ the more general case can then be obtain upon replacing $\delta$ by $\delta^\alpha$.
We  consider the reverse H\"older inequality
$$ \|g\|_q \|\delta^s\|_p<\|g \delta^s\|_1$$
for $0<p<1$ and $q=\frac{p}{p-1}$. 
%Suppose that $s<\underline{\eta}\alpha^{-1}$ and 
Let $p=\frac1{\underline{\eta}-\epsilon}$, 
$\epsilon>0$.
Then ${g^\frac{p}{p-1}}
={g^\frac{-1}{\underline{\eta}-\epsilon-1}} $  is summable and hence 
$\|g\|_q>0$.
Thus $\|g \delta^s\|_1<+\infty$ implies
$\|\delta^s\|_{p}<+\infty$. But the latter is
equivalent to $\frac{s}{\underline{\eta}-\epsilon}>1$.
It follows that
 $s_0^{low}\geq \underline{\eta}$. 

Given that $(\delta_n)_n\in\ell^{1+\epsilon}$ there exists
$c>0$ such that $\delta_n^{s}\leq c
n^{\frac{-s}{1+\epsilon}}=\frac{c}{n^{{s}-\epsilon'}}$, $\epsilon'=\epsilon'(\epsilon)>0$.
Hence
$$g_n\delta_n^s\leq
\frac{g_n}{n^{s-1-2\epsilon'}}\frac{c}{n^{1+\epsilon'}},$$
%$\epsilon'=\epsilon'(\epsilon)>0$ such that 
%$\epsilon'\stackrel{\epsilon\to 0}{\to} 0$.
Note that $c/ n^{1+\epsilon'}$ is summable. Hence
if $s>\overline{\eta}+2\epsilon'$ then
$\frac{g_n}{n^{s-1-2\epsilon'}}$ is 
bounded and, by the standard H\"older inequality with $p=1$,
$\|g \delta^s\|_1<+\infty$. It follows that $s_0^{low}\leq
\overline{\eta}$. 
\qed

\medskip
The above theorem formulates bounds on $s_0^{low}$ in terms of
properties of the sequence $g$. This
sequence does not take into account the spacing between $r_n$ and
$r_{n+1}$. While this does not matter in the case of subshifts for
instance,
for tilings where the spacing between $r_n$ and
$r_{n+1}$ can vary greatly, 
it is more customary to define $\delta_n$ by a function
$\delta$ over $\R^+$ through the formula  $\delta_n=\delta(r_n)$.
We will therefore formulate a second result which puts restriction on
the integrability of the function $\delta$. We start with a theorem from \cite{PB09}.
For that we anticipate Lemma~\ref{lem-PB-zeta} of the last section stating that the abscissa of 
convergence of the zeta function $\zeta^{PB}$ of the spectral triples of Pearson-Bellissard coincides with $s_0^{low}$ in case that the Michon graph has uniformly bounded branching.
%%%%%%%%
\begin{thm}[\cite{PB09}]
Suppose that the Michon graph associated to a compact ultra metric space $(X,d)$
has uniformly bounded branching. Then the abscissa of 
convergence $ s_0^{PB}$ of the zeta function $\zeta^{PB}$ is the upper box dimension:
$$ s_0^{PB} = \overline{\dim}(X,d) = \overline{\lim}_n \frac{\ln \#\Tt^{(0)}_n}{-\ln \delta_n}.$$
Here $\#\Tt^{(0)}_n$ denotes the number of vertices of level $n$ in the Michon tree.
\end{thm}
%%%%%%%%
Note that we can reformulate $(\delta_n)_n\in \ell^s(\s)$ as
$\delta\in L^s(dP)$ where $\delta$ is any  
continuous function satisfying $\delta(r_n)=\delta_n$ and $dP$ the
(atomic) measure on $\R$ having $P$ as cumulant, i.e.\ $dP$ is the
sum over $n$ of $\s(n)$ times the Dirac measure at $r_n$. 

The following is relatively easy to prove, given that for a tiling $\#\Tt^{(0)}_n$ is equal to the value of the patch counting function at $r_n$.
%%%%%%%%
\begin{lemma}[\cite{Ju09}] Consider a $d$-dimensional tiling of finite local complexity and provide its discrete tiling space with the metric defined by $\delta(r) = \frac{1}{r}$.
Then
%%%
\[
\overline{\beta} = \overline{\dim}(\Xi,d)\quad\mbox{and}\quad
\underline{\beta} = \underline{\dim}(\Xi,d)\,,
\]
%%%
where $\underline{\dim}$ is the lower box dimension
 \end{lemma}
%%%%%%%%
This allows to conclude as follows
%%%%%%%%
\begin{cor}
Consider %the function $\zeta^{low}$ of 
a $d$-dimensional tiling of finite local complexity and with a uniform bound on the number of possible patch extensions.
Let $\delta$ be a piecewise $C^1$-function with negative derivative defining the metric on $\Xi$.
Suppose that $\delta\in L^{1+\epsilon}([1,+\infty))$ (w.r.t.\ Lebesgue measure) for all $\epsilon>0$.
Then $$s_0\leq\overline{\beta}.$$
\end{cor}
%%%%%%%%
{\em Proof:} $\delta\in L^{1+\epsilon}([1,+\infty))$ implies $r\delta(r)\leq c$ for some $c>0$. 
Hence, $\ln r \leq -\ln\delta(r) + \ln c$ implying   $s_0=\overline{\dim}(\Xi,d)=\overline{\beta} \,\overline{\lim}_n \frac{\ln r_n}{-\ln\delta(r_n)}\leq \overline{\beta} $.\qed

\medskip
A relation between the abscissa of convergence and the complexity exponent was first made in \cite{JS10a} where it is proved that $s_0=\beta$ if $\delta(r) = \frac{1}{r}$ in the context of
(primitive) substitution tilings of $\RM^d$. In this context one
can prove further that $s_0$ is the Hausdorff dimension of the discrete tiling space \cite{JS10c}.

In the following theorem we obtain also a lower bound on $s_0$ even without the assumption that the tiling has a uniform bound on the number of possible patch extensions.
%%%%%%%%%%%%%%
\begin{thm}
\label{thm-zeta-2}
Consider a $d$-dimensional tiling of finite local complexity.
Let $\delta$ be a piecewise $C^1$-function with negative derivative defining the metric on $\Xi$.
Suppose that
$\delta\in L^{\alpha+\epsilon}([1,+\infty))\backslash
L^{\alpha-\epsilon}([1,+\infty))$ (w.r.t.\ Lebesgue measure) for some $\alpha >0$ and for all $\epsilon>0$ (small enough).
Then
$$\alpha\underline{\beta} \leq s_0^{low} \leq \alpha\overline{\beta}.$$
Moreover  
$$\alpha \underline{\beta} \leq s_0 \leq \alpha(\overline{\beta}+d-1)$$
where $s_0$ is the abscissa of convergence of the zeta-function of the tiling.
\end{thm}  
%%%%%%%%%%%%%%
{\em Proof:} 
Again, we may concentrate on the case $\alpha=1$.
We have, up to irrelevant additive constants arising from
the boundary terms in the partial integration,
$$\zeta^{low}(s) = \int_{r_0}^\infty \delta^s dP =  - \int_{r_0}^\infty
\big(\delta^s(r)\big)' P(r) dr .$$
Since $\delta'\leq 0$ this yields, for all $\epsilon>0$ 
$$- \int_{r_0}^\infty
\big(\delta^s(r)\big)' 
r^{\underline{\beta}-\epsilon} dr \leq \zeta^{low}(s) \leq 
- \int_{r_0}^\infty
\big(\delta^s(r)\big)' 
r^{\overline{\beta}+\epsilon} dr.$$ 
After partial integration this yields, again  up to irrelevant
additive constants, 
$$\int_{r_0}^\infty \delta^s(r) r^{\underline{\beta}-1-\epsilon} dr \leq \zeta^{low}(s)\leq 
\int_{r_0}^\infty \delta^s(r) r^{\overline{\beta}-1+\epsilon} dr.$$ 
In particular $$ \sup\{s:\int_{r_0}^\infty \delta^s(r)
r^{\underline{\beta}-1-\epsilon} dr=+\infty\}\leq s_0^{low}\leq
\inf\{s:\int_{r_0}^\infty \delta^s(r)
r^{\overline{\beta}-1+\epsilon} dr<+\infty\}.$$
We thus need to study the integrability of $ \delta^s(r)
r^{\gamma}$. As before, $\delta\in L^{1+\epsilon}([1,+\infty))$
implies that $\delta\leq {c}{r^{\frac{-1}{1+\epsilon}}}$. Hence
$ \delta^s(r) r^{\gamma}\leq
{c}{r^{\gamma+\frac{-s}{1+\epsilon}}}$ and the latter is
integrable if $\frac{s}{1+\epsilon}>\gamma+1$. We apply this to 
$\gamma = \overline{\beta}-1$ to find that
$s_0^{low}\leq  \overline{\beta}$. 

To study the values of $s$ 
which imply the non-integrability of $ \delta^s(r)
r^{\gamma}$ we use again the reverse H\"older inequality as in the proof of Theorem~\ref{thm-zeta-1}. 
The role of
$g$ is played by the function $g(r) = r^{\gamma}$. With that $g$ the
corresponding $\underline{\eta}$ must be $\gamma+1$. We apply this to
$\gamma =\underline{\beta}-1$. 
The arguments parallel to those in the last proof yield 
$s_0^{low}\geq  \underline{\beta}$.

To get the upper bound for $s_0$ we bound $a(v)+1$ by $r_n^{d-1}$ 
where $n$ is the level of $v$ ($r_n^{d-1}$ is the
surface of the boundary of an $r$-patch). We thus need to
study the summability of $ \sum_n g(n) r_n^{d-1}
\delta_n^{s}$. This shifts $\overline{\beta}$ by $d-1$.
\qed

\medskip
For tilings, the two pairs of exponents $\beta$ and $\eta$ are in general not easily related, and this comes from the lack of control over the interspacings between the $r_n$'s.
For subshifts however the exponents are related as follows.
%%%%%%%%%%%%%%%%%%
\begin{prop}
\label{prop-zeta-2}
For a subshift with $(\delta_n)_n \in \ell^{\alpha+\epsilon}\backslash \ell^{\alpha-\epsilon}$ for some $\alpha >0$ and for all $\epsilon>0$ (small enough), one has
%%%%%%%%%%%%%%%%%%
\[
\underline{\beta} \le \underline{\eta} \le s_0^{low}/\alpha \le \overline{\eta} \le \overline{\beta}\,. 
\]
%%%%%%%%%%%%%%%%%%
\end{prop}
%%%%%%%%%%%%%%%%%%
{\em Proof:}
By Theorem~\ref{thm-zeta-1} it suffices to show \(\underline{\beta} \le \underline{\eta}\) and \(\overline{\eta} \le \overline{\beta}\).
For the first inequality, note that \( g^{-1/(\gamma - 1)}(n)\) is summable for all \(\gamma>\underline{\eta}\), and since it is decreasing, one has \( g(n) \ge c_\gamma n^{\gamma -1}\) for some constant $c_\gamma>0$.
Now \(P(n)-P(1) = \sum_{j=1}^n g(j) \ge c_\gamma \sum_{j=1}^n n^{\gamma-1}\), and this last sum is bounded below by the integral \(\int_{1}^n r^{\gamma-1} dr \ge c n^\gamma\).
So we get \(P(n) \ge c n^\gamma\) for all \(\gamma>\underline{\eta}\) and the first inequality follows.

For the second, one has \( g(n) \le c_\gamma n^{\gamma-1}\) for all \(\gamma \le \overline{\eta}\), and we similarly get \(P(n) -P(1) \le c_\gamma \sum_{j=1}^n n^{\gamma-1}\) which is bounded above by the integral \(\int_1^n r^{\gamma-1} dr \le c n^\gamma\).
\qed

%\newpage
%%%%%%%%%%%%%%%%%%%%%%%%%%%%%%%%%%%%%%%%%%%%%%%%%%%
\subsection{Laplacians}

We aim to define Laplacians which have domains which intersect densely with 
%are defined on (dense subspaces of) $L^2$-functions 
continuous functions on $X$. 
The graph Laplacians will not do it as they are defined on
$\ell^2$-functions on $V$ and no continuous 
non-zero function will restrict to an $\ell^2$-function on $V$.
The way to solve this problem is to average over the choices. For that we
consider a probability space $(Y,\P)$ for the choice functions. 
Recall from Section~\ref{approxgraph} that a choice function $\tau : \Tt^{(0)} \rightarrow \Pi_\infty$ is uniquely determined by the
family of functions $v\mapsto \tau_v(v):=(\tau(v))_{|v|+1}$.
We may therefore identify $Y$ with the Cartesian product  
$Y = \Pi_{v\in \Tt^{(0)}} Y_v$ where $Y_v$ is the (finite) set of choices $\tau_v:v\to \Tt^{(0)}(v)$. 
%%%%%%%%%%%%%%%%
\begin{lemma}
\label{lem-measure}
There is a one-to-one correspondence between Borel probability measures on
$X$ and  product probability measures
$d\P(\tau) = \Pi_{v\in \Tt^{(0)}}d\P_v(\tau_v)$ on $Y$. Indeed, given a family of
probability measures $\P_v$ on $Y_v$ 
\[
\mu([v]) := 
%\EE(\chi_{v\:*}) =
\Pi_{n=0}^{|v|-1}\P_{v_n}\{\tau_{v_n}(v_n)=v_{n+1}\}
\] 
defines a probability measure on $X$, where $v_n$ is the n-th vertex on the (unique) path from the root $\circ$ up to $v$.
Conversely, given a Borel probability measure $\mu$ on $X$ 
\[
\P_{v}\{\tau_{v}(v)=u\} := \frac{\mu([u])}{\mu([v])}
\]
yields a probability measure on $Y_v$.
\end{lemma}
%%%%%%%%%%%%%%%%
The proof is straightforward. Equivalently we may say that
$$ \int_X f(x) d\mu(x) = \int_Y f(\tau(\circ))d\P(\tau) $$
for all continuous functions on $X$.

From now on we will therefore suppose that our measure $\P$ is a product measure
and denote by $\mu$ its corresponding measure on $X$.
This allows us to define the (real)
Hilbert space $L^2(X,\mu)$. 
We denote its scalar product $ \langle f,g \rangle$.
%%%%%%%%
\begin{lemma}
For any $n \in \NM$
$$\langle f,g \rangle =   \sum_{v\in \Tt^{(0)}_n} 
%\EE((\chi_v){\circ}) 
\mu([v])\int_Y f(\tau(v)) g(\tau(v)) d\P(\tau). $$ 
%Furthermore, $ \EE(\chi_v) = \Pi_{n=0}^{|v|-1}\P_{v_n}(v_{n+1})$.
\end{lemma}
%%%%%%%%
{\em Proof:} 
Using Lemma~\ref{lem-measure} we have
%%%%%%%%%%%%%%%%%%
\[
\langle f, g \rangle = \int_X f(x) g(x) \; dx = \int_Y f(\tau(\circ)) g(\tau(\circ)) \; d\PM(\tau)\,.
\]
%%%%%%%%%%%%%%%%%%
We can write
\(f(\tau(\circ))  = \sum_{u\in \Tt^{(0)}_n} f(\tau(u)) \chi_u (\tau(\circ))\), where $\chi_u$ is the characteristic function of the set of infinite paths going through $u$.
We get
%%%%%%%%%%%%%%%%%%
\[
\langle f, g \rangle = \sum_{u,v\in \Tt^{(0)}_n}  \int_Y f(\tau(u)) g(\tau(v))  \, \chi_u (\tau(\circ)) \chi_v(\tau(\circ)) \;d\PM(\tau)\,.
\]
%%%%%%%%%%%%%%%%%%
Now \(\chi_u (\tau(\circ)) \chi_v(\tau(\circ)) = \chi_v(\tau(\circ))\) if $u=v$, and $0$ otherwise, so we deduce
%%%%%%%%%%%%%%%%%%
\[
\langle f, g \rangle = \sum_{v \in \Tt^{(0)}_n} \int_Y f( \tau(v) ) g(\tau(v)) \, \chi_v(\tau(\circ)) \; d\PM(\tau)\,.
\]
%%%%%%%%%%%%%%%%%%
We now remark that \(\chi_v(\tau(\circ)) = \Pi_{n=0}^{|v|-1} \chi_{v_{n+1}}(\tau_{v_n}(v_n))\), where $v_n, n \le |v|$, is the (unique) vertex at level $n$ on the path from the root up to $v$.
And we use the product decomposition of $Y$ as 
\(Y(v) \times \Pi_{n=0}^{|v|-1} Y_{v_n} \times Y^{(v)}\), 
where $Y(v)$ is the product of $Y_u$'s over all vertices $u$ in the subtree starting at $v$, and $Y^{(v)}$ is the product of the remaining $Y_u$'s (for all $u\neq v_0 =\circ, v_1, v_2 \cdots v$ and not in the subtree starting at $v$).
We get
%%%%%%%%%%%%%%%%%%%
%\begin{multline*}
%\int_Y f( \tau(v) ) g(\tau(v)) \, \chi_v(\tau(\circ)) \; d\PM(\tau) = 
%\int_{Y(v)} f( \tau(v) ) g(\tau(v)) \; d\PM(v) \\
%\Pi_{n=0}^{|v|-1} \int_{ Y_{v_n}} \chi_{v_{n+1}}(\tau_{v_n}(v_n))  d\PM_{v_n}
%\int_{Y^{(v)}} d\PM^{(v)}
%\end{multline*}
%%%%%%%%%%%%%%%%%%%
%%%%%%%%%%%%%%%%%%
\begin{eqnarray*}
\langle f, g \rangle & = & \sum_{v \in \Tt^{(0)}_n} \int_{Y(v)} f( \tau(v) ) g(\tau(v)) \; d\PM(v) \quad \cdots \\
& & \qquad \qquad \cdots \quad
%\int_{\Pi_{n=0}^{|v|-1} Y_{v_n}} \chi_{v}(\tau(\circ))  d\PM \; 
\Pi_{n=0}^{|v|-1} \int_{ Y_{v_n}} \chi_{v_{n+1}}(\tau_{v_n}(v_n))  d\PM_{v_n} \; \int_{Y^{(v)}} d\PM^{(v)} \\
%& & \qquad \qquad \qquad \int_{Y(v)} f( \tau(v) ) g(\tau(v)) \; d\PM(v) \\
& = & \sum_{v \in \Tt^{(0)}_n} \int_{Y(v)} f( \tau_v(v) ) g(\tau_v(v)) \; d\PM(v)  \quad \Pi_{n=0}^{|v|-1}\P_{v_n}\{\tau_{v_n}(v_n)=v_{n+1}\} \\
& = &   \sum_{v \in \Tt^{(0)}_n} \mu([v]) \int_{Y} f( \tau(v) ) g(\tau(v)) \; d\PM(\tau) \,,
\end{eqnarray*}
%%%%%%%%%%%%%%%%%%
where we used Lemma~\ref{lem-measure} to identify $\mu([v])$.
\qed
\medskip

We wish to define a quadratic form on a dense subspace of $L^2(X)$
by the formula
%%%
\[
Q_\rho(f,g) = \frac{1}{2}\int_Y \Tr_{\ell^2(E)}
(\rho(D) [D,\pi_\tau(f)]^*[D,\pi_\tau(g)])d\P(\tau)\,,
\]
%%%
where $\rho$ is a positive function (which we see as a ``density matrix'').
It defines a closable Dirichlet form \cite{FOT94} whose domain is generated by (real valued) locally constant functions on $X$ (the reader can find a proof of this in a similar setting in the work of Pearson-Bellissard \cite{PB09}).
Since $Q_\rho$ is positive and symmetric this gives rise to
a densely defined symmetric operator $\Delta_\rho$ by the formula
%%%
\[
 Q_\rho(f,g) = \langle \Delta_\rho f,g\rangle 
\]
%%%
which we wish to determine.

First note that 
$$ Q_\rho(f,g) = \sum_n g_n  q_n(f,g)$$ 
where $g_n = \rho(\delta_n)\delta_n^{-2}$,
%$$ q_n(f,g) = \int_Y \Tr_{\ell^2(E_n)}
%([D_n,\pi_\tau(f)]^*[D_n,\pi_\tau(g)])d\P(\tau).$$  
\[
q_n(f,g) = \frac{1}{2} \sum_{(u,v)\in H^{(1)}_n} \int_Y 
(f_{u}-  f_{v})({g_{u}-  g_{v}}) d\P,
\]
% where $r(e)_n,s(e)_n\in \Tt^{(0)}_n$ denote the representative vertices of  $r(e), s(e) \in V$, and for $u \in \Tt^{(0)}$ 
where we have used the notation  $f_u:Y\to\R$, $f_u(\tau) = f(\tau(u))$ for any $u \in \Tt^{(0)}$.
In fact, for $n$ fixed any function on $X$ gives rise to a finite
family of functions $(f_u)_{u\in \Tt^{(0)}_n}$ which have the property that
$f_u(\tau) = f_u(\tau')$ for all $\tau'$ such that $\tau'(u)=\tau(u)$.
Conversely, any family of functions $(f_u)_{u\in \Tt^{(0)}_n}$ with that
property determines a function on $X$. 

We will obtain  $\Delta_\rho$ as a sum 
$\Delta_\rho = \sum_n g_n \Delta_n$ 
with  $  q_n(f,g) = \langle \Delta_n f,g\rangle $.
%%%%%%%%%%
\begin{lemma}
\label{lem-Delta} 
For all $f\in L^2(X)$
%%%%%%%%%%%%%%%%%%%
\[
\Delta_n f(x) = \frac{1}{\mu([x_n])}
\sum_{y : (x,y)\in E_n}(f(x)-\EE(f_{y_n}))\,,
\]  
%%%%%%%%%%%%%%%%%%%
where $x_n, y_n \in \Tt^{(0)}_n$ denote the vertices of $x,y$ at level $n$ under the identification $X\cong \Pi_\infty$.
\end{lemma}
%%%%%%%%%%
{\em Proof:} Note that $(x,y) \in E_n$ implies $(x_n,y_n) \in H^{(1)}_n$ so that we have
%%%%%%%%%%%%%%%%%%%
\[
q_n(f,g) = \sum_{v\in \Tt^{(0)}_n} \mu([v]) \int_Y
\sum_{u : (u,v) \in H^{(1)}_n} \mu([v])^{-1}(f_v(\tau) - f_u(\tau)) g_v(\tau)d\P(\tau).
\]
%%%%%%%%%%%%%%%%%%%
The function $F_{v,u} = f_v - \EE(f_u)$ satisfies that
$F_{v,u}(\tau')=F_{v,u}(\tau)$ for $\tau'$ such that
$\tau'(v)=\tau(v)$. Hence
\((\Delta_nf)_v=\sum_{u:(u,v)\in H^{(1)}_n}\mu([v])^{-1} F_{v,u}\) provided that 
$ \int_Y f_u(\tau) g_v(\tau)d\P(\tau) =\int_Y \EE(f_u)
g_v(\tau)d\P(\tau) $ for all $u\neq v$. 
To show this let decompose $Y =Y(v)\times Y'$ where $Y(v)$ is the product of the $Y_{u'}$'s for all vertices $u'$ is the subtree starting at $v$, and $Y'$ the product of the remaining $Y_{u'}$'s.
Likewise let $d\PM(\tau) = d\PM(v) d\PM'$. 
Then
\begin{eqnarray*} 
\int_Y \EE(f_u) g_v(\tau)d\P(\tau) & =& 
\int_{Y(v)\times Y'} \EE(f_u) g_v \;d\PM(v) d\PM'\\
& = & \int_{Y'} \EE(f_u) \;d\PM' \ \int_{Y(v)} g_v \;d\PM(v) \\
& = & \EE(f_u)\EE(g_v)\,,
\end{eqnarray*}
and similarly, writing $Y= Y(v) \times Y(u) \times Y''$, one gets
\begin{eqnarray*} 
\int_Y f_u(\tau) g_v(\tau) d\P(\tau) & =& 
\int_{Y(v)\times Y(u) \times Y''} f_u g_v \; d\P(v) d\P(u) d\P''\\
& =& 
\int_{Y(v)} g_v \; d\P(v) \ \int_{Y(u)}  f_u \; d\P(u) \ \int_{Y''} d\P''\\
& = & \EE(f_u)\EE(g_v) \,,
\end{eqnarray*}
from which the claim  follows. \qed
\bigskip

Given a vertex $u \in \Tt^{(0)}$ we denote by $\chi_u$ the characteristic function of the set of infinite  paths in $\Pi_\infty \cong X$ that go through $u$.
For a finite rooted path $\gamma$, we denote by $\chi_\gamma$ the characteristic function of the set of infinite paths with prefix $\gamma$ (that is the set of paths through the last vertex of $\gamma$).
For $n\in\NM$ and an infinite path $x$ (resp. a finite path $\gamma$ of length $|\gamma| \ge n$) we let $x_n$ (resp. $\gamma_n$) denote the vertex at level  $n$  through which $x$ (resp. $\gamma$) goes. 

It is useful to express the action of $\Delta_\rho$ on characteristic functions $\chi_\gamma$. 
%%%%%%%%%%
\begin{thm} 
We have
 $\Delta_\rho 1 = 0$ and
%%%
\[
\Delta_\rho\chi_\gamma =
\sum_{n=1}^{|\gamma|}\frac{\rho(\delta_n)}{\mu([\gamma_n])\delta_n^{2}}\left( 
a(\gamma_{n-1})\chi_\gamma
- \mu([\gamma]) \sum_{u:(u,\gamma_n)\in H^{(1)}_n} \frac{1}{\mu([u])}\chi_{u}\right)\,,
\]
%%%
where \(a(\gamma_{n-1})=\# \Tt^{(0)}(\gamma_{n-1}) -1 \) is the branching number of $\gamma_{n-1}$ minus 1. 
\end{thm}
%%%%%%%%%%
{\em Proof:} 
We apply Lemma~\ref{lem-Delta} to $\chi_\gamma$.
Using that $\chi_\gamma(\tau(u)) =1 $ whenever $\tau(u) = \tau(\gamma)$
we find
%%%%%%%%%%%%%%
\[
\EE((\chi_{\gamma})_u) = \left\{\begin{array}{ll}
1 & \mbox{if $\gamma$ is a prefix of the path from $\circ$ up to $u$} \\
p_n(\gamma) & \mbox{if } u = \gamma_n, \quad n\leq|\gamma| \\   
0 & \mbox{otherwise} 
\end{array}\right.
\]
%%%%%%%%%%%%%%
where \(p_n(\gamma) = \PM\{ \tau(\gamma_n) = \tau(\gamma)\} = \mu([\gamma])/\mu([\gamma_n])\) by Lemma~\ref{lem-measure}.
First let $n>|\gamma|$ and $(u,v)\in H^{(1)}_n$.
Then $\chi_\gamma(\tau(v))=1$ iff $\gamma$ is a prefix of the path from $\circ$ up to $v$. % the path from $\circ$ to $v$ is a extension of $\gamma$ iff 
It follows that $\chi_\gamma(\tau(v))=\EE((\chi_{\gamma})_u)$ and hence $\Delta_n(\chi_\gamma) = 0$.

Now let $n\leq|\gamma|$ and $(u,v)\in H^{(1)}_n$.
If $v=\gamma_n$ then $u\neq\gamma_n$ and $\EE((\chi_{\gamma})_u)=0$.
If $v\neq\gamma_n$ but  $v \in \Tt^{(0)}( \gamma_{n-1})$ then one and only one of the $u$'s satisfies $u=\gamma_n$.
Hence in the last case $\sum_{u:(u,v)\in H^{(1)}_n} \EE((\chi_{\gamma})_u)=p_n(\gamma)$.
In all other cases none of the  $u$'s satisfies $u=\gamma_n$.
It follows that
\begin{eqnarray*}
\frac1{\mu([x_n])} \sum_{u:(u,x_n)\in H^{(1)}_n}\EE((\chi_{\gamma})_u)
&=&\frac{p_n(\gamma)}{\mu([x_n])}(\chi_{\gamma_{n-1}}(x_n)-\chi_{\gamma_n}(x_n))\\
&=&\frac{\mu(\gamma)}{\mu([\gamma_n])}\sum_{u:(u,\gamma_n)\in H^{(1)}_n}
\frac1{\mu([u])} \chi_u(x_n).
\end{eqnarray*}
\qed

%%%%%%%%%%%%%%%%%%%%%%%%%%%%%%%%%%%%%%%%%
\subsection{Comparison with the construction of Pearson-Bellissard}
\label{sec-comparison} 
There is some flexibility in the above construction and this gives us
the possibility to compare our construction with the one given by
Pearson and Bellissard in \cite{PB09}. 

Let $\E$ be a function on the set of branching vertices
$\Tt^{(0)}_{br} \subset \Tt^{(0)}$ that associates to $v$ a selection
of pairs of distinct elements in $\Tt^{(0)}(v)$, that is, a symmetric subset of 
$\Tt^{(0)}(v)\times \Tt^{(0)}(v)\backslash\diag(\Tt^{(0)}(v)) $.
Viewing again a pair of vertices as an edge we define a new
neighborhood graph $H(\E) = (\Tt^{(0)}, H^{(1)}(\E))$ which has the
  same
vertices as the original one of section~\ref{sec-3.1} and edges
 $H^{(1)}_n(\E) = \bigcup_{v\in\Tt^{(0)}_{n-1}} \E(v)$. Using the same
 equivalence relation as in that section we thus obtain a graph
$\Gamma(\tau,\E)$  with the same vertices as $\Gamma(\tau)$ but
possibly fewer edges.
Clearly  $\Gamma(\tau,\E)=\Gamma(\tau)$ if $\E(v)$ consists of all
possible pairs. 
On the other hand, if $\E(v)$ contains a single unoriented pair (that is, two oppositely oriented pairs) then the
spectral triple associated to $\Gamma(\tau,\E)$ coincides with one of
Pearson-Bellissard's. 
To explain this we describe their triple in our framework.

Define a PB-choice function to be a function
$$\tau_{PB}=(\tau^+,\tau^-):\Tt^{(0)}_{br}\to X\times X\backslash\diag(X)$$
such that for $v\in \Tt^{(0)}_{n-1}$, we have $\tau^\pm(v)\in [v]$ and the pair $(\tau^+(v)_n, \tau^-(v)_n) $ forms a horizontal edge in $H^{(1)}_n$. 
A PB-choice function defines a representation $\pi_{\tau_{PB}}$ of
$C(X)$ on $\ell^2(\Tt^{(0)}_{br})\otimes\C^2$, 
%%%
\[
\pi_{\tau_{PB}} (f)\Psi(v) = \left(\begin{array}{cc}f(\tau^+(v)) & 0 \\
0 & f(\tau^-(v))\end{array}\right)\Psi(v)\,, \quad \forall v \in \Tt^{(0)}_{br}\,.
\]
%%%
Furthermore, define a Dirac operator by
\[ D_{PB}\Psi(v) = \frac1{\delta_{n-1}}\left(\begin{array}{cc}
0 & 1 \\ 1 & 0 \end{array}\right) \Psi(v), \quad v\in  \Tt^{(0)}_{n-1}\cap \Tt^{(0)}_{br}\, . \]
The family of spectral triples obtained in this way is the family of Pearson-Bellissard spectral triples.
We compare them with
the spectral triple we associated to $\Gamma(\tau,\E)$ where
$$
\E(v) = \{(\tau^+(v)_n, \tau^-(v)_n),(\tau^-(v)_n, \tau^+(v)_n) \}\,, \quad v\in \Tt^{(0)}_{n-1} \cap \Tt^{(0)}_{br}\,,
$$ 
and $\tau$ is a choice function of our type which ought to satisfy
\begin{equation}\label{eq-choice}
\tau(\tau^+(v)_{|v|+1}) = \tau^+(v),\quad \mbox{and}\quad 
\tau(\tau^-(v)_{|v|+1})= \tau^-(v).
\end{equation}  
Moreover, the PB-choice function induces a decomposition of the edges $E(\E)$ of $\Gamma(\tau,\E)$, namely $E(\E)^+$ contains precisely the edges of the form $(\tau^+(v)_n, \tau^-(v)_n)$ and $E(\E)^-$ those of the form $(\tau^-(v)_n, \tau^+(v)_n)$.

While equations (\ref{eq-choice}) cannot be guaranteed for all vertices, as our choice
functions satisfy additional restrictions, the converse holds: given
a choice function of our type, a function $\E$ assigning to each
branching vertex $v$ a single pair of oppositely oriented edges,
and a choice of orientation for the pair in $\E(v)$,
one can easily construct a PB-choice function so that the above equations hold. 
 Consider in this case the unitary 
\begin{eqnarray*}
U\:\::  \quad
%\ell^2(E(\E)^+)\oplus \ell^2(E(\E)^-)
\ell^2(E(\E))\quad &\to &\ell^2(\Tt^{(0)}_{br})\otimes\C^2 \\
1_{(\tau^+(v)_n, \tau^-(v)_n)} & \mapsto & 1_v \otimes \left(1 \atop 0 \right) \\
1_{(\tau^-(v)_n, \tau^+(v)_n)} & \mapsto & 1_v \otimes \left(0 \atop 1 \right) 
\end{eqnarray*}
Then $U^{-1} \pi_{\tau_{PB}}(f) U = \pi_\tau (f)$ and $U^{-1} D_{PB} U = D$.
Also $U$ preserves the $\Z_2$-grading, i.e.\ the spectral triples are unitary equivalent.
It follows in particular that when the Michon tree has the property
that each vertex branches into at most two vertices then our spectral
triples are particular members of the family of PB-spectral triples. 

In the case that vertices may branch into more than two vertices the
graph  $\Gamma(\tau,\E)$ will no longer be connected, moreover the
connected components of the graph have closures in $X$ which do not
intersect and hence the spectral distance associated with Pearson-Bellissard spectral
triple does not even yield a metric (it takes value $+\infty$). 
This is why they do have to take all choices into
account and, for instance, take the infimum over all choice functions
as in Corollary~\ref{cor-inf} to obtain a metric.  
Strictly speaking, therefore, the spectral distance associated to any
single one of their spectral triples is of limited use.
\bigskip

Denote by $\zeta^{PB}(s) = \Tr(|D_{PB}|^{-s})$ the zeta function of Pearson Bellissard
associated to the spectral triple (we ignore the factor $\frac12$ they put reflecting the obvious double
degeneration of the spectrum of $D_{PB}$
due to the tensor product with the Hilbert space with $\C^2$).   
%%%%%%%%%%%%
\begin{lemma}
\label{lem-PB-zeta}
$\zeta^{PB}\leq \zeta^{low}$. If the Michon tree has uniformly
bounded branching
then $\zeta^{PB}$, $\zeta^{low}$, and $\zeta$ have all the same
abscissa of convergence.
\end{lemma}
%%%%%%%%%%%%
{\em Proof:} As in \cite{PB09} one gets
$$\frac12\zeta^{PB}(s) = \sum_n \sum_{v\in \Tt^{(0)}_{n}\cap
  \Tt^{(0)}_{br}} \delta_n^s = \sum_n\sum_{v\in \Tt^{(0)}_{n}}
\bar{b}(v) 
\delta_n^s $$ 
where $\bar{b}(v)=1$ if $v$ branches and $0$ otherwise.
In particular  $\bar{b}(v)\leq a(v)$ from which the first claim
follows. By definition, the Michon tree has uniformly bounded
branching if
 $a(v)\leq c \bar{b}(v)$ for some constant $c\in\N$. Hence the
 second claim.
\qed
\bigskip

Given that the Dirichlet forms are in both approaches defined via an
average over choice functions, it is not surprising that they look 
a lot alike.
To  compare them we describe the PB-Dirichlet form. 
%%%%
\[
Q_\rho^{PB} = \sum_n \frac{\rho(\delta_n)}{\delta_n^2} \sum_{v\in
 \Tt^{(0)}_{n-1}} q_v^{PB}
\] 
%%%%
where
%%%%
\[
q_v^{PB}(f,g) = \int_{\Ypsilon}
\Big(f(\tau^+(v))-f(\tau^-(v))\Big) \Big(g(\tau^+(v)-g(\tau^-(v))\Big)
d\P^{PB}(\tau)
\]
%%%%
where $(\Ypsilon,\P^{PB})$ is the probability space of PB-choice
functions. 
$\P^{PB}$ is a product measure on  $\Ypsilon=\Pi_{v\in \Tt^{(0)}_{br}} \Ypsilon_v$.
It is uniquely determined by a measure $\mu$ on $X$ through the formula
%%%%
\[
\int_{\Ypsilon_v} F(\tau^+(v),\tau^-(v)) d\P^{PB}_v(\tau) =
N(v)^{-1}\!\!\!\!\!\!\! 
\sum_{(u_1,u_2)\in \E(v)} \int_{[u_1]\times[u_2]} \!\!\!\!   F(x,y)d\mu(x)d\mu(y)
\]
%%%%
for all continuous functions $F$ on $X\times X$. The normalization
constant is given by 
$N(v) = \sum_{(u_1,u_2)\in \E(v)}\mu([u_1])\mu([u_2])$.

Since $N(v) = 2\mu([u_1])\mu([u_2])$ provided $u_1,u_2$ are the only branching vertices of $v$ we find exactly the same expression as for our Laplacian in the case that any branching vertex splits into only two vertices. 

In the case where there are more than two branching vertices, however,
we can obtain the PB-Laplacian only if we include 
an average over the choices for the $\E(v)$. 
Let $\nu_v$ be a symmetric probability measure on the set $\E(v)$.
Then $\tilde Q_\rho = \sum_{n}\frac{\rho(\delta_n)}{\delta_n^2}\sum_{v\in \Tt^{(0)}_{n-1}}\tilde q_v$ where 
%%%%%
\[
\tilde q_v(f,g) = \frac{1}{2}\sum_{(u_1,u_2)\in \E(v)} \nu_v((u_1,u_2))\int_Y
(f_{u_1}- f_{u_2})(g_{u_1}- g_{u_2})d\P
\]
%%%%%
reproduces the PB-Laplacian, provided \(\nu_v((u_1,u_2))= \mu([u_1])\mu([u_2])N(v)^{-1}\).

%Likewise we may consider the zeta-function associated to $\Gamma(\tau,\E)$ which depends on $\E$ but not on $\tau$.
%If we require that $\E(v)$ is not empty (and symmetric) then $a(v)\geq 1$ for branching vertices and therefore $\xi^{low}\leq \xi_{\Gamma(\tau,\E)}\leq \xi$.
%Furthermore, $\xi^{low} = \xi_{\Gamma(\tau,\E)}$ if $\E(v)$ consists of one un-oriented edge. 

\end{document}